\documentclass[ijoc,nonblindrev]{informs3b}

\OneAndAHalfSpacedXI
\usepackage{natbib}
 \bibpunct[, ]{(}{)}{,}{a}{}{,}%
 \def\bibfont{\small}%

\TheoremsNumberedThrough
\ECRepeatTheorems
\EquationsNumberedThrough
\usepackage{amssymb,amsmath,mathtools,algorithm,graphicx, color, gensymb, xfrac, stmaryrd, soul, multirow, epstopdf}
\usepackage{subcaption}
\usepackage[T1]{fontenc}
\usepackage[dvipsnames]{xcolor}
\usepackage{graphicx}
\usepackage{pbsi,varwidth}
\usepackage[countmax]{subfloat}
\usepackage[flushleft]{threeparttable}

\usepackage[mathscr]{euscript}
\usepackage[left, mathlines, displaymath]{lineno}
\usepackage[noend]{algpseudocode}

\usepackage{url}
\usepackage{nicefrac}
\usepackage{amssymb}
\usepackage{amsmath}
\usepackage{hyperref}
\hypersetup{
	colorlinks,
	linkcolor={red!50!black},
	citecolor={blue!50!black},
	urlcolor={blue!80!black}
}
\usepackage{breakurl}
\newcommand{\LeftEqNo}{\let\veqno\@@leqno}
\usepackage{environ}
\usepackage{xcolor}
\usepackage{float}
\usepackage{graphicx}
\usepackage[noend]{algpseudocode}
\usepackage{algorithm}
\usepackage{tikz}
\usepackage{float}
\usepackage{array}
\usepackage{enumerate}
\usepackage{subcaption}
\usepackage{booktabs}
\usepackage{comment}
\usepackage{stackrel}
\usepackage{color,soul}
\usepackage{algorithmicx}
\usepackage{graphicx}
\usepackage{xspace}
\usepackage{lscape}
\usepackage{tikz}
\usepackage{bbm}

\usepackage[margin=1in]{geometry} % required 1in margin at four sides

\newcommand*\patchAmsMathEnvironmentForLineno[1]{%
	\expandafter\let\csname old#1\expandafter\endcsname\csname #1\endcsname
	\expandafter\let\csname oldend#1\expandafter\endcsname\csname end#1\endcsname
	\renewenvironment{#1}%
	{\linenomath\csname old#1\endcsname}%
	{\csname oldend#1\endcsname\endlinenomath}}%
\newcommand*\patchBothAmsMathEnvironmentsForLineno[1]{%
	\patchAmsMathEnvironmentForLineno{#1}%
	\patchAmsMathEnvironmentForLineno{#1*}}%
\AtBeginDocument{%
	\patchBothAmsMathEnvironmentsForLineno{equation}%
	\patchBothAmsMathEnvironmentsForLineno{align}%
	\patchBothAmsMathEnvironmentsForLineno{flalign}%
	\patchBothAmsMathEnvironmentsForLineno{alignat}%
	\patchBothAmsMathEnvironmentsForLineno{gather}%
	\patchBothAmsMathEnvironmentsForLineno{multline}%
	\patchBothAmsMathEnvironmentsForLineno{eqnarray}%
}

\usepackage{natbib}
\setcitestyle{authoryear,open={(},close={)}}

\newcommand{\APP}{{AP-SRMU}}
\newcommand{\PRM}{{SRMU}}
\newcommand{\PLUB}{{PLUB}}
\newcommand{\MIP}{{MIP}}

\begin{document}

%	\linenumbers
	\RUNAUTHOR{Mehrani, S., Sefair, J.A.}
%	\RUNAUTHOR{-----------------}
	\RUNTITLE{Robust Assortment Optimization under Sequential Ranking-based Choice Model}
	\TITLE{Robust Assortment Optimization under a Ranking-based Choice Model  with Product Unavailability Effect}
	\ARTICLEAUTHORS{%
		\AUTHOR{Saharnaz Mehrani}
		\AFF{Department of Operations and Information Management, School of Business, University of Connecticut, Storrs, CT 06269, \EMAIL{saharnaz.mehrani@uconn.edu}}
		\AUTHOR{Jorge A. Sefair}
		\AFF{School of Computing, Informatics, and Decision Systems Engineering, Arizona State University, Tempe, AZ 85281, \EMAIL{jorge.sefair@asu.edu}} %, \URL{}}
	}
	
	\ABSTRACT{
	The assortment planning problem is a central piece in the revenue management strategy of any company in the retail industry. In this paper, we study a robust assortment optimization problem for substitutable products under a sequential ranking-based choice model. Our modeling approach incorporates the cumulative effect of finding multiple unavailable products on the customers' purchase decisions. To model the highly uncertain order in which a customer explores the products to buy, we present a bi-level optimization approach to maximize the expected revenue under the worst-case order of products in the preference lists of customers. We provide a polynomial-time algorithm that optimally solves a special case of the unconstrained assortment planning problem under our choice model. For the general constrained version of the problem, we devise a solution procedure that includes a single-level reformulation and a cutting-plane approach to iteratively tighten the solution space. We also provide a greedy algorithm that can quickly solve large instances with small optimality gaps.}

\maketitle 
\section{Introduction}\label{sec:Introduction}
The assortment planning problem is a central piece in the revenue management strategy of every retail company. In this problem, retailers decide which products to offer in order to match the customers' preferences while maximizing the total revenue. Offering a small product assortment may result in a loss of potential sales, customer dissatisfaction, and a loss of goodwill. On the contrary, offering a large set of products, aiming to capture all possible demands, requires significant investments in space, inventory, and logistics. To decide the optimal assortment, retailers must identify the critical aspects of the customers' purchasing behavior and then translate them into a prescriptive decision model, which are typically based on mathematical programming due to the combinatorial nature of the assortment planning problem \citep{kok2008assortment, besbes2016product}.

%In this work, we focus on modeling some aspects of the customers' purchasing behavior that have been overlooked in the literature. We model the assortment decisions as a function of the inherent value of each product for the customers, the expected revenue from product sales, correlation among the demands for different products, and the customers' sensitivity to product unavailability.
Features like price and quality are commonly used factors that determine the value of each product for customers and influence product demand. However, if a product is not offered in the assortment, customers may substitute that product with an alternative available product, an effect known as \textit{substitution behavior} \citep{blanchet2016markov}. As a result, the demand of an offered product not only depends on its own features, but also on the whole set of products in the assortment and the strength of the customers' substitution behavior. 

The substitution behavior is determined by the product's attributes, customers' behavior, and the cost of searching for the preferred products at other stores. According to \citet{anderson2006measuring}, substitution is negligible for durable products like bedding and home accessories in the moderate-to-premium price range, which implies that customers can be very sensitive to product unavailability. Moreover, customers negatively react to product unavailability when the personal commitment to a product increases and when it is possible to buy the product from other sellers (either at a brick-and-mortar store or online) \citep{fitzsimons2000consumer}. The substitution effect is especially strong nowadays given the increasing penetration of online retailing and the resulting low cost of finding preferred products. As a result, customers may be more sensitive to the unavailability of their preferred products. However, for some types of products such as apparel and furniture, with high variety in style, color, and quality, the search cost for inspecting different stores is very high, encouraging customers to make in-store substitutions to other products in the assortment. In some cases, product unavailability in a store can even facilitate the customer's purchase decision \citep{wang2017impact, ma2019assortment}.

We propose a robust assortment optimization problem under a sequential ranking-based model that captures \textit{the cumulative effect of product unavailability}, where the probability of leaving the store with no purchase increases with the number of substitution attempts due to product unavailability. We use the general ranking-based model, as it has the ability to capture this effect in the sequence of products examined by the customer. The general ranking-based model assumes that each customer has a preference list of products, whose pre-established order dictates the customer purchasing decisions (see \citet{honhon2012optimal} and \citet{goyal2016near} for details). The order of the products in the list is based on the utility maximization mechanism, i.e., the list contains products in descending order of utility and all products in the list provide utilities higher than the no-purchase option  \citep{mahajan2001stocking}. We use the same construction of the preference list, but allow customers to have different sensitivity levels to the unavailability of their top-choice products. Depending on the product type and the customer's tolerance to unavailability, a customer may or may not make a substitution when attempting to purchase a product whose similar or preferred alternatives are easy to find in other stores. As a result, customers may leave the store without making any purchase after finding that some of their preferred products are unavailable. Throughout the paper, we refer to the probability that the customer leaves the store with no purchase as the \textit{leaving probability}.

\cite{fitzsimons2000consumer} shows in a laboratory experiment that consumers respond more aggressively to unavailability of products of higher preference, resulting in a high likelihood of switching store. Therefore, when customers substitute all their top priority products, it means they are determined to make a purchase as long as they can obtain a utility higher than making no purchase. Our model also considers those cases where the customer is determined to make a purchase in the current store, even after a number of substitutions representing the end of the top-priority product list. The length of such list depends on the product category. To this end, we partition each customer's preference list into two levels, top and low priority, and only consider the cumulative effect when top-priority products are unavailable. 

We present an assortment optimization problem to maximize the expected revenue generated by products in the assortment, while also considering the uncertainty in the choice model. The uncertainty is primarily driven by the utility of products for each customer, which determines their position in the preference list. Because it is not easy to infer such preference lists, we adopt a \textit{robust} approach that focuses on the worst-case list that maximizes the leaving probability of a given assortment. This list consists of the worst-possible order of products from the retailer's perspective. However, we control the level of conservatism in this robust approach by only considering the worst-case list for top-priority products. We use \textit{Multinomial Logit} (MNL) utilities for low-priority products, which is a tractable and the most common random utility model.

We formulate the robust approach using a bi-level optimization model that selects an assortment that maximizes the retailer's expected revenue (first level) while facing customers with a priority list of maximum leaving probability (second level). We develop exact and greedy optimization algorithms to solve the assortment planning problem under our sequential %\textcolor{red}{partial} 
ranking-based choice model. %\textcolor{red}{In particular, we provide a polynomial-time algorithm for the unconstrained assortment problem when the worst-case customer's preference list contains only one unavailable product.}
In particular, we provide a polynomial-time algorithm when the assortment problem is unconstrained and there is only one top-priority product in every preference list. Our solution procedure for the general unconstrained problem and that under a cardinality constraint includes a transformation into a single-level upper-bound problem, which we improve with a cutting-plane procedure. We also devise a set of problem-specific valid inequalities that can improve the performance of our solution algorithms. Additionally, we develop a greedy (heuristic) algorithm that can solve large-scale problem instances efficiently in terms of running time and whose empirical testing shows promising results in terms of the achieved optimality gap. Indeed, we use the results from this (fast) greedy algorithm as an initial (feasible integer) solution to improve the performance of our exact solution methods. Moreover, we extend our models and algorithms to multiple customer categories with different choice parameters. In this case, we show that ignoring the variability in the sensitivity to product unavailability across customer segments can result in low-revenue assortments. %\textcolor{red}{We also present numerical evidence to characterize the sensitivity of the optimal assortment to changes in the length of the customers' preference lists for different types of products.}
We also present numerical evidence to characterize the sensitivity of the optimal assortment to changes in the number of top-priority products in a preference list.

The remainder of this paper is organized as follows. Section \ref{sec:LiteratureReview} provides a review on the related literature. We describe our sequential
ranking-based choice model under product unavailability in Section \ref{sec:ChoiceModel} and our bi-level assortment optimization approach in Section \ref{sec:Model}. Section \ref{sec:Solution} presents our solution algorithms. We numerically show the effectiveness of our solution methods and analyze the performance of the proposed choice model in Section \ref{sec:Results}. We present our final remarks and future work in Section \ref{sec:Conclusion}.

\section{Literature review}\label{sec:LiteratureReview}

In recent decades, assortment optimization problems under different choice models have received considerable attention in the operations research literature. Existing works focus on finding a balance between the generality of the proposed customer choice models and the tractability of the resulting optimization problems. For this reason, this section presents a literature review consisting of two parts: choice models and solution methods.

\subsection{Choice models}\label{sec:LiteratureReview-choice}

In particular, our work is related to the assortment planning problem under both the general ranking-based choice model and its special case, MNL. The use of ranking-based choice models in assortment planning was first presented by \citet{mahajan2001stocking} and later extended by \citet{rusmevichientong2006nonparametric} for multi-product pricing. \citet{aouad2018approximability} study an assortment planning problem under a general ranking-based choice model, where customer choices are modeled through an arbitrary distribution over a set of ordered product preference lists. The authors argue that it is difficult to accurately approximate general ranking preferences and provide best-possible approximability bounds. \citet{jagabathula2016nonparametric} examine a ranking-based choice model in a non-parametric joint assortment planning and pricing problem where customers first select a subset of products with prices below a predetermined threshold and then choose the most preferred product out of the resulting set. \citet{smith2000management} study the impact of having customers' preference lists of different lengths on the demand distribution, which they relate to the maximum number of substitution attempts. The authors show that assortments with a large number of products are less sensitive to the increase in the number of substitution attempts, which otherwise deteriorate the assortment performance. Moreover, allowing multiple substitution attempts makes this assortment planning problem more difficult to solve, justifying the development of specialized solution methods for the simplified single-substitution attempt case \citep{kok2007demand}.

Because the ranking-based choice model is more general than many others, \citet{honhon2010assortment} and \citet{goyal2016near} use this choice model in joint assortment planning and inventory management problems. \citet{honhon2010assortment} consider different types of customers, each of which has a specific ordered preference list of products. They assume that the proportion of each type of customers is known and that the retailer can estimate the possible preference lists and their associated probabilities. \citet{goyal2016near} study a joint assortment and inventory planning problem under dynamic (stock-out-based) substitution with a simple consumer choice model, which is later extended to a customer choice model that depends on both price and quality. The authors also assume that the number of arriving customers follows an increasing failure rate distribution and that the customers' preference lists have different sizes with a known probability distribution.

The assortment planning problem for substitutable products under different variants of the MNL choice model has received considerable attention in the literature (see e.g., \citet{ryzin1999relationship}, \citet{talluri2004revenue}, \citet{rusmevichientong2010dynamic}, \citet{kok2015assortment}). \citet{rusmevichientong2010dynamic} study static and dynamic assortment planning problems under the MNL choice model and subject to a cardinality constraint that limits the maximum number of products chosen. In the static version of the problem, they assume that the parameters of the MNL model are known in advance, while in the dynamic problem these are estimated from data. %They propose a polynomial time algorithm to solve the static version of the problem and develop a policy to estimate the unknown parameters and simultaneously solve the dynamic problem. 
More recent extensions of the MNL model, such as the mixed MNL \citep{bernstein2015dynamic,feldman2015bounding, kunnumkal2015upper,kunnumkal2019tractable} and the nested logit model \citep{gallego2014constrained, davis2014assortment}, incorporate more complex customer's choice behaviors. \citet{bernstein2015dynamic} study a problem in which a customized assortment of products is offered to customers that arrive sequentially. In this case, products are identically priced and substitutable with limited inventory in a market with multiple customer segments. They use the mixed MNL choice model in which each customer belongs to a segment with a probability and each segment follows a specific MNL model. %Aiming to maximize the expected revenue under the mixed MNL, \citet{feldman2015bounding} present a solution approach that provides tight upper bounds on the total revenue, considering the problem complexity.
\citet{gallego2014constrained} study an assortment optimization problem with cardinality and space constraints under the nested logit choice model in which products are categorized in different nests. In such choice model, products in each nest are substitutable and customers first select a nest, and then a product within the selected nest \citep{davis2014assortment}. \citet{flores2019assortment} solve an assortment planing problem under the sequential MNL in which products are partitioned into two levels such that customers only examine the products in the second level if they do not purchase a product in the first level. For a sample of the extensive literature on the MNL choice model, see \citet{ryzin1999relationship}, \citet{kok2015assortment}, \citet{feldman2015bounding}, \citet{kunnumkal2019tractable}, and \citet{liu2019assortment}.

In an effort to generalize the choice models, \citet{blanchet2016markov} present a Markov chain-based approximation for all random utility based discrete choice models, including MNL, probit, nested logit, and mixed MNL. In this approximation, they define a Markov chain with states given by the candidate products (as well as the no-purchase option) and whose transition matrix captures the substitution probabilities. These probabilities only depend on the last visited product due to the Markovian memoryless property, which implies that the cumulative effect of finding multiple products unavailable is limited to one product. \citet{berbeglia2016discrete} generalizes the Markov chain-based model by proposing a random walk-based model that solves the memoryless limitation by considering that the substitution probability depends on the whole sequence of previously visited products.  %\citet{desir2019constrained} use the choice model from \citet{blanchet2016markov} and solve an assortment optimization problem under the Markov chain-based choice model and subject to capacity constraint by developing constant factor approximation solution algorithms.  

\subsection{Solution methods}\label{sec:LiteratureReview-Opt}

Multiple exact and heuristic approaches have been proposed for the solution of different variants of the assortment planning problem. \citet{aouad2018approximability} study an unconstrained assortment optimization under the general ranking-based choice model. The authors relate the problem to the maximum independent set problem and prove that it is NP-hard to approximate within a $O(n^{1-\epsilon})$ factor for any $\epsilon > 0$ and $n$ candidate products. They present best-possible approximation algorithms with performance guarantees as a function of both extremal prices and the maximum length of any preference list. \citet{honhon2010assortment} and \citet{goyal2016near} study joint assortment and inventory planning problems under the general ranking-based model and stochastic demand. \citet{honhon2010assortment} focus on the unconstrained problem and present a dynamic programming (DP) algorithm to find multiple local maxima using the properties of their DP value function. %The quality of their heuristic algorithm improves as the proportion of customers over different preference lists is random. 
Their experimental results on 4, 6, and 8 candidate products and different distributions and parameter values show no optimality gap in most cases under random preference lists  and an improvement in the run time compared to existing heuristics from the literature. \citet{goyal2016near} provide NP-hardness results of the capacity-constrained problem even for the special case where there is only one customer and all preference lists consist of only two products. Under some assumptions, including that each customer has a price threshold and always prefers the cheapest available product, they develop a polynomial-time approximation scheme (PTAS) for computing near-optimal solutions with arbitrary level of accuracy. %They also show how to leverage this scheme for the case that the choice model depends on both price and quality. 
Their experimental results on up to 20 candidate products and different capacity limits show that PTAS is able to produce assortments with at least $44\%$ of the optimal expected revenue in 7200 seconds.

The lack of complete data and the difficulty to estimate the input parameters of some assortment planning models have motivated the use of bi-level and robust optimization  approaches to guarantee a maximum revenue even under worst-case scenarios of customer behavior. Bi-level optimization has a rich literature in modeling competitive markets and has been recently used to formulate problems in the area of revenue management in the airline industry \citep{birbil2009role, perakis2010robust} and portfolio optimization  \citep{chen2018robust}. In assortment optimization, \citet{rusmevichientong2012robust} develop a bi-level formulation of the assortment optimization problem under the MNL model in order to guarantee a high revenue under unknown choice model parameters. The authors consider both static and dynamic settings of the problem, where there can be a limited initial product inventory that must be allocated over time, and prove revenue-ordered characterizations of the optimal assortment in both settings. \citet{rooderkerk2016robust} develop a robust optimization framework to protect the assortment against the demand and profit uncertainty for every product, achieving a balance between risk and return. The authors study a capacity-constrained problem  and provide a heuristic that approximates the risk-return efficient frontier of assortments. \citet{li2019robust} present a robust assortment optimization problem under the MNL choice model considering partially available information for the parameters of the choice model. The assortment model uses a chance constraint that requires the revenue to exceed a given threshold value with a fixed probability, which is then replaced by a distributionally robust chance constraint. This chance constraint is also approximated using the worst-case Conditional Value-at-Risk. \citet{desir2019nonconvex} use a bi-level optimization formulation for an unconstrained assortment optimization problem under the general Markov chain-based choice model and its special case, MNL, considering the uncertainty in the Markov chain's transition probabilities. They maximize the worst-case expected revenue, where the worst-case is taken over a set of likely values of the choice model parameters. The authors show that under certain assumptions related to the choice model and the uncertainty set, the order of the max and min operators can be interchanged in the bi-level objective function. Using this result, the authors develop an iterative algorithm that converges in polynomial time under few conditions.

This paper contributes to the existing literature in four main aspects. First, although our choice model can be seen as a special case of the more general Markov-chain and random-walk based models of \citet{blanchet2016markov} and \citet{berbeglia2016discrete}, respectively, and the general ranking-based models of \citet{mahajan2001stocking} and \citet{farias2013nonparametric}, our approach explicitly incorporates the \textit{cumulative} effect that the number of unavailable products attempted (i.e., number of substitution attempts) has on the customers' leaving probability. %This allows us to devise a mechanism to calculate the worst-case leaving probability as a function of the products excluded from the assortment. Moreover, we transform the ranking-based choice model into a two-step process, where each customer first follows an ordered preference list and then, if the customer reaches the end of the ordered list and decides not to leave the store, examines the offered assortment following the MNL choice model.
Moreover, the sequential process in our choice model allows us to include those cases where the customer is determined to make a purchase in the current store. Second, we propose a bi-level optimization approach which is designed to find an optimal assortment that provides the maximum possible revenue under the worst-case customers' preference list and their tolerance to unavailable products. We call this strategy ``robust'' as it is intended to produce an optimal assortment with limited and highly uncertain information on the customers' preferences. Third, we present both exact and heuristic solution approaches to solve multiple variants of the resulting bi-level assortment planning problem. These approaches are designed to both overcome the difficulty of solving the nonlinear integer bi-level problem formulation and take advantage of the problem structure. We present a polynomial-time algorithm for a special case of the unconstrained assortment problem. To solve the general constrained versions of the problem, our algorithmic approach first transforms the bi-level problem into a single-level upper bound problem. Using problem-specific properties, we devise a cutting-plane approach and problem-specific valid inequalities that iteratively tighten the single-level upper bound formulation. Additionally, we use the results from the unconstrained case to construct super-valid equalities to further accelerate the solution of special cases of the constrained problem. Our approaches can solve large-scale instances (e.g., $>80$ candidate products) of some variants of the problem, which is larger than the instances currently solved in the literature. Fourth, we extend our models and algorithms to multiple customer categories with different choice parameters. In this context, we show that ignoring the variability in the sensitivity to unavailable products across customer segments can result in low-revenue assortments. We also present numerical evidence to characterize the sensitivity of the optimal assortment to changes in the length of the customers' preference lists for different types of products.

\section{Sequential ranking-based choice model under product unavailability effect} \label{sec:ChoiceModel}
%\textcolor{red}{In this section, we formulate a \emph{partial ranking-based choice model that considers the effect of product unavailability} on the customers' purchasing behavior, which we refer to as \PRM{}.}
In this section, we present our \emph{sequential ranking-based choice model that considers the effect of product unavailability} on the customers' purchasing behavior, which we refer to as \PRM{}. We define the set of all 
%\textcolor{red}{candidate} 
products 
%\textcolor{red}{for the assortment}
by $ I=\{0, 1,...,n\} $, where $0$ denotes the no-purchase option. We assume that 
%\textcolor{red}{all}
products in $I$ belong to the same category and that they are substitutable with each other. Moreover, we assume that each customer purchases at most one product in this category. In our model, each customer entering the store has a preference list of products ranked in non-increasing order of utility. 
%\textcolor{red}{In our model, each customer entering the store has a preference list of products ranked in decreasing order of preference}. \textcolor{red}{We d}
Denote such list by $ \mathcal{T} = \{\tau_1,...,\tau_l\} \subseteq I $, where $ \tau_1 \succeq ... \succeq \tau_l $  ($ a \succeq b $ denotes that $ a $ is at least as preferred as is $ b $). If the utility of option $ i \in I $ for the customer with preference list $\mathcal{T}$ is $ U_i $, then $ U_{\tau_1} \geq ... \geq U_{\tau_l} \geq U_0 $, where $U_0$ is the utility assigned to the no-purchase choice, regardless of the assortment. %\textcolor{red}{When product $ \tau_k $ is unavailable, for any $ k \in \{1,...,l-1\} $, the customer either chooses the no-purchase option and leaves the store or attempts to substitute $ \tau_{k} $ with $ \tau_{k+1} $. If product $ \tau_{l} $ is unavailable, then the customer either leaves the store or purchases a product from the assortment according to the MNL.}

Because of the customer sensitivity to the unavailability of the products with the highest utility, we partition the preference list of each customer into two priority levels. The top-priority products are the top $\bar u$ products in the preference list, while the remaining products are defined as low-priority. If product $ \tau_k $ is unavailable, for any $ k \in \{1,...,\bar u\} $, i.e., any top-priority product in a preference list $\mathcal{T}$, then customer either attempts to substitute $ \tau_{k} $ with $ \tau_{k+1} $ or chooses the no-purchase option from the current retailer. This means that the customer does not necessarily make a substitution when examining the top-priority list. When the customer finds all the products in the top-priority list unavailable and remains in the store (i.e., after $\bar{u}$-1 substitutions), it means that the customer wants to make a purchase and therefore, buys the product with the largest utility among the available ones.

We model the substitution probability of top-priority product $\tau_k$ as a function of two factors: product features and the position in the preference list. Let $ p_i^0 $ be the leaving probability (i.e., the probability of choosing the no-purchase option) when $ i $ is the first product in the preference list, i.e., $ \tau_1=i $, and this product is unavailable. To capture the customer's dissatisfaction due to product unavailability, we assume that the leaving probability increases as the customer makes more substitution attempts, moving down the preference list. Mathematically, we define $ \eta_{ik} $ as the rate of increase in the leaving probability $ p_i^0 $, when the customer has faced $ k $ unavailable top-priority products and the last attempted purchase was product $ i $. Accordingly, the customer substitutes unavailable product $ \tau_k=i $ with product $ \tau_{k+1} $ with probability $ 1-\eta_{ik}p_i^0 $. We define $ \eta_{i1}=1 $ for all $ i \in I $, meaning that the probability of substituting the first product in the preference list, $ \tau_1=i $, is $ 1-p_i^0 $. Parameter $\eta_{ik}$ is useful to model the sensitivity of customers to unavailable products. For instance, a rapidly increasing value of $\eta_{ik}$ with respect to $k$ indicates that either customers are very sensitive to the unavailability of their top choices and will not attempt many substitutions or that products in the given category are readily available in other stores (i.e., low search cost). %\textcolor{blue}{Therefore, we also assume that parameters $ p_i^0 $ and $\eta_{ik}$ are exogenous and independent from product features and utility $U_i$ for any $i \in I$.}

Define the set of products in the assortment by $ S \subseteq I \setminus \{0\}$ and the set of unavailable products by $ \bar S= I \setminus (S \cup \{0\})$. Denote the set of all possible preference lists by $L$ and the probability of observing a customer with preference list $\mathcal{T} \in L$ by $w_{\mathcal{T}}$ with $\sum_{\mathcal{T} \in L} w_{\mathcal{T}} = 1$. Given an assortment $S$, the probability of purchasing product $i \in S$ is given by

\begin{equation}
   \boldsymbol{p}(i,S) = \sum_{\mathcal{T} \in L} w_{\mathcal{T}} \times
    \begin{cases}
        1 & \text{if $i \in {\mathcal{T}}^1, \tau_1=i$} \\
        \prod_{k=1}^{k(i)-1} (1-\eta_{\tau_k k} \ p_{\tau_k}^0) & \text{if $i \in {\mathcal{T}}^1, \{\tau_1,...,\tau_{k(i)-1}\} \subseteq \bar S$} \\
        \prod_{k=1}^{\bar u} (1-\eta_{\tau_k k} \ p_{\tau_k}^0)  & \text{if $i \in {\mathcal{T}}^2, \{\tau_1,...,\tau_{k(i)-1}\} \subseteq \bar S$} \\
        0 & \text{otherwise,} \label{e:purchprob}
    \end{cases}
\end{equation}
where $k(i)$ is the position of product $i$ in the preference list $\mathcal{T}$ and $\mathcal{T}^1$ and $\mathcal{T}^2$ denote the sets of top- and low-priority products, respectively. The term $1-\eta_{\tau_k k} \ p_{\tau_k}^0$ is the substitution probability for the unavailable product in position $k$ of preference list $\mathcal{T}^1$. The first condition in \eqref{e:purchprob} means that product $i$ is purchased with certainty if it is the first in the list ${\mathcal{T}}$. Otherwise, it is purchased as long as all the more preferred products in the list ${\mathcal{T}}$ are unavailable and the customer decides to substitute them. In this case, the purchase probability is $\prod_{k=1}^{k(i)-1} (1-\eta_{\tau_k k} \ p_{\tau_k}^0)$ if $i \in {\mathcal{T}}^1$ and   $\prod_{k=1}^{\bar u} (1-\eta_{\tau_k k} \ p_{\tau_k}^0)$ if $i \in {\mathcal{T}}^2$. Example \ref{E:ChoiceModel} illustrates the construction of the purchase probability.

\begin{example} \label{E:ChoiceModel} Let $ n=2 $, $L=\{ \{1\},\{2\},\{1,2\},\{2,1\} \}$ with respective probability vector $w=\{ 0.1,0.2,0.3,0.4 \}$, $ p_1 ^0=0.3 $, $ p_2^0=0.2 $, and $\bar u =1$. If $ S= \{1\} $, the probability of purchasing Product $1$, $\boldsymbol{p}(1,S)$, is $ w_1+w_3+w_4 \ (1-p_{2}^0) = 0.72$. The no-purchase probability, $\boldsymbol{p}(0,S)$, is $1-0.72=0.28$, which can also be calculated as $w_2+w_4 \ p_{2}^0$. 
\end{example}

We mitigate the uncertainty on the actual distribution of the utility vector $U$ by proposing a robust approach that focuses on the worst possible preference list of products, which is the list that maximizes the leaving probability. This also helps us eliminate the dependency of \eqref{e:purchprob} to the  distribution of $L$, which is challenging to obtain from empirical data. To be less conservative in our optimization, we find the worst possible preference list only for top-priority products. The details of this construction are explained in Section \ref{sec:Model}. Given a worst-case top-priority list, we describe the purchasing decisions of low-priority products using the MNL model, which is a widely used special and tractable case of the general ranking-based model \citep{mahajan2001stocking}.

%\textcolor{red}{When the customer attempts to buy all the products in the preference list and finds all of them unavailable, then either chooses the no-purchase option or a product with the maximum utility among those in the assortment. That is, we assume that if all products in the preference list are unavailable and the customer decides to make more substitution attempts, then the customer behaves according to the MNL. Under the MNL model, each customer assigns a utility $ U_i $ to each option $ i \in I $, which consists of a deterministic component, $ u_i $, and random component, $ \varepsilon_i $, that follows a Gumbel distribution with location and scale parameters $\mu$ and $\beta$, respectively. That is, $ U_i=u_i+\varepsilon_i $, $ \forall i \in I $ \citep{kok2008assortment}. We assume that the $ \varepsilon $-random variables are independent and identically distributed with $\mu=0$ and $\beta=1$, as in \citet{rusmevichientong2010dynamic}. We also assume that $ u_0=0 $, which means that the no-purchase option provides zero expected utility to the customers. Defining the set of products in the assortment by $ S \subseteq I \setminus \{0\}$, the customer chooses the product with the highest positive utility among those in $ S $. As a result, the probability that a customer chooses option $ i \in S \cup \{0\} $ under the MNL choice model is given by} 

In the MNL, $U_0$, $U_1$, ...,$U_n$ are mutually independent random variables with Gumbel distributions with means $u_0$, $u_1$, ...,$u_n$, respectively \citep{kok2008assortment}. Similar to \citet{rusmevichientong2010dynamic}, we assume that the scale parameter of all distributions is equal to 1. We also assume that $ u_0=0 $, which means that the no-purchase option provides zero expected utility to the customers. Under the MNL, the probability that a customer chooses product $ i $ given assortment $S$ is given by
\begin{equation}\label{equ:equ00}
\rho_i(S)=\frac{\nu_i}{\nu_0+\sum_{j \in S} \nu_j},
\end{equation}
where $ \nu_i=e^{u_i} $ for all $ i \in S \cup \{0\} $ \citep{kok2008assortment}. If $\hat{\mathcal{T}}^1=\{\tau_1,...,\tau_{\bar u}\}$ is the worst possible preference list that customers can have for their top-priority products, then the probability of purchasing product $i \in S$ is given by
\begin{equation} 
   \boldsymbol{p}(i,S) = 
    \begin{cases}
        1 & \text{if $i \in \hat{\mathcal{T}}^1, \tau_1=i$} \\
        \prod_{k=1}^{k(i)-1} (1-\eta_{\tau_k k} \ p_{\tau_k}^0) & \text{if $i \in \hat{\mathcal{T}}^1, \{\tau_1,...,\tau_{k(i)-1}\} \subseteq \bar S$} \\
        \prod_{k=1}^{\bar u} (1-\eta_{\tau_k k} \ p_{\tau_k}^0) \ \rho_i(S)  & \text{if $i \notin \hat{\mathcal{T}}^1, \{\tau_1,...,\tau_{\bar u}\} \subseteq \bar S$} \\
        0 & \text{otherwise.} \label{e:Robustpurchprob}
    \end{cases}
\end{equation}

To obtain \eqref{e:Robustpurchprob} from  \eqref{e:purchprob}, we replace $w_{\mathcal{T}}$ by using the MNL probabilities for low-priority products. That is because the top-priority list (i.e., $\hat{\mathcal{T}}^1$) is known, as it corresponds to the worst-case sequence (i.e., $w_{\mathcal{T}}$ is replaced with 1 for the first two conditions). As a result, we only need to find the probability of purchasing the low-priority product $i \in S$, which corresponds to $\rho_i(S)$. % and  which replaces $w_{\mathcal{T}}$ for the third condition in Equation \ref{e:purchprob}.}

\section{Bi-level assortment optimization} \label{sec:Model}
In this section, we develop a robust mathematical formulation for the \emph{assortment planning problem under the \PRM{}} from Section \ref{sec:ChoiceModel}, which we refer to as \APP{}. We formulate the \APP{} as a bi-level optimization problem that maximizes the expected revenue under the worst-case customer preference list. The motivation of this formulation is the difficulty in finding the preference list probability distribution among different customers interested in buying the same type of product \citep{aouad2018approximability}. 

We define the worst-case customer preference list for a given assortment of products as a sequence of unavailable products that minimizes the probability of staying in the store (or maximizes the probability of leaving the store without any purchase). We control the length of the preference list using parameter $\bar u$, which is the maximum number of allowed top-priority products in the preference list. A large value of $\bar u$ indicates that customers have many product options to explore at the store in case of unavailability of more preferred items. If $\bar u = 0$ our model reduces to the MNL because the customer immediately faces the decision of selecting a product from the assortment or leave the store. We assume that the value of $\bar{u}$ is less than the number of products not included in the assortment so it is always possible to find a preference list consisting of $\bar{u}$ unavailable products. To find the worst-case top-priority list of length $\bar u$, we define the binary decision variable $ y_{ik} $, which is equal to $1$ if product $ i $ is located in position $ k $ in the preference list, and is equal to $0$ otherwise, for $i \in I \setminus \{0\}$ and $\ k \in \{1,...,{\bar u} \} $. Given an assortment $ S $, the set of feasible worst-case preference lists is defined by
\begin{equation}\label{equ:equ03}
\begin{aligned}
Y(S)= \left\{ \mathbf{y} \in \{0,1\} ^ {{\lvert \bar S \rvert} \times {\bar u}} :\ \sum_{i \in \bar S}y_{ik}= 1, \ \forall k \in \{1,...,{\bar u}\}, \ \sum_{k=1}^{\bar u}y_{ik}\leqslant 1, \ \forall i \in \bar S \right\},
\end{aligned}
\end{equation}
where $ \bar S= I \setminus (S \cup \{0\})$ as we only consider unavailable products to be part of a worst-case preference list. The worst-case preference list is the one that minimizes the probability of staying in the store while attempting to purchase unavailable products, which is given by

\begin{equation}\label{equ:equ3}
\begin{aligned}
\pi(S)=\min_{\mathbf{y} \in Y(S)} \prod_{k=1}^{{\bar u}}\left(1-\sum_{i \in \bar S} {\eta_{ik} p_i^0 y_{ik}}\right).
\end{aligned}
\end{equation}

We assume that if $ S = I\setminus\{0\} $ (i.e., $ \bar S = \emptyset $) in Problem (\ref{equ:equ3}), then every customer can purchase their top choice and therefore $ \pi(I \setminus \{0\}) = 1 $. We also assume that $ \pi(S)=1 $ when $ \bar u=0 $. Intuitively, having less products in the assortment increases the chances that a customer attempts to buy a product that is unavailable. This is formalized in the following lemma, where we show that removing products from any assortment $ S \subseteq I \setminus \{0\} $ never increases the value of $ \pi(S) $.
\begin{lemma} \label{l:lambda00}
	$ \pi(S^ \prime) \leq \pi(S) $ for any $ S \subseteq I \setminus \{0\} $ and $ S^\prime \subseteq S $.
\end{lemma}

\proof{Proof.} 
Define $ \bar S^\prime = I \setminus (S^\prime \cup \{0\}) $ and $\bar S = I \setminus (S \cup \{0\})$, thus $ \bar S \subseteq \bar S^\prime $. Denote the optimal solution to Problem (\ref{equ:equ3}) for assortment $S$ by $ y^\ast_{ik}(S) $ for $ i \in \bar S $ and $ k \in \{1,...,\bar u\} $. Construct a feasible solution to Problem (\ref{equ:equ3}) for assortment $S^\prime$ as follows: $ y^\prime_{ik}(S^\prime) = y^{*}_{ik}(S) $ for $ i \in \bar S $ and $ k \in \{1,...,\bar u\} $ and $ y_{ik} = 0 $ for all $ i \in \bar S^\prime \setminus \bar S $ and $ k \in \{1,...,\bar u\} $. Therefore, we have
\begin{align*}
\pi(S^ \prime) & \leq \prod_{k=1}^{{\bar u}}\left(1-\sum_{i \in \bar S^\prime} {\eta_{ik} p_i^0 y^\prime_{ik}(S^\prime)}\right) \\
& \leq \prod_{k=1}^{{\bar u}}\left(1-\sum_{i \in \bar S} {\eta_{ik} p_i^0 y^\ast_{ik}(S)}\right) = \pi(S),
\end{align*}
where the second inequality holds because $ \bar S \subseteq \bar S^\prime $ and $\eta_{ik} p_i^0 \geq 0$ for all $i \in I$ and $k \in \{1,...,\bar u\}$, proving the statement. $ \square  $ 

We define $ r_i $ as the revenue from product $ i \in I $ and let $ \varOmega $ be the set of feasible assortments, which can include different types of constraints. To maximize the expected revenue considering the worst-case customers' preference list under the \PRM{}, the retailer must solve the optimization problem
\begin{equation}\label{equ:equ4}
\begin{aligned}
z^\ast= & \max _{S \in \varOmega} \left\{\pi(S) \sum_{i \in S} \rho_i(S)r_i\right\}. 
\end{aligned}
\end{equation}
Note that the \APP{} in (\ref{equ:equ4}) gives the retailer the possibility to capture different customer behaviors by using specific values of $ \bar u $ and the $\eta$-parameters. 

\section{Solution approach} \label{sec:Solution}
In this section, we describe a procedure to find an optimal solution to the \APP{} in $ (\ref{equ:equ4}) $. In Section $ \ref{subsec:UnconSolution} $, we present a polynomial-time solution approach for the unconstrained version of the problem (i.e., $\varOmega = \{S: S \subseteq I \setminus \{0\}\}$) when $ \bar{u}=1 $. In Section $ \ref{subsec:ConSolution} $, we present a reformulation and a solution strategy for Problem $ (\ref{equ:equ4})$ under a cardinality constraint that limits the number of products in the assortment. In the same section, we provide an alternative greedy algorithm for solving large-scale instances. We also describe a strategy to strengthen the cardinality-constrained version of Problem $ (\ref{equ:equ4})$ when $\bar u = 1$ based on the solution of the unconstrained problem. At the end of this section we discuss how our solution approaches can be used to solve some extensions of the \APP{}.

\subsection{Polynomial algorithm for unconstrained \APP{} and $\bar{u} = 1$} \label{subsec:UnconSolution}   
In this section, we focus on the unconstrained version of Problem $ (\ref{equ:equ4}) $ (i.e., $ \varOmega= \{S:S \subseteq I \setminus \{0\} \} $) and where $ \bar u = 1 $. This can be the case for many durable products such as furniture and some electronic devices for which customers have a single preferred product in mind. Nevertheless, due the high cost of searching for alternatives at other stores, some customers may decide to explore alternative products from the assortment in case of unavailability of the only product in their preference list. The unconstrained assumption reflects the case of stores selling products on demand (or by catalog), where customers see store samples before placing an order for product shipment. This strategy allows the store to have a very large (virtually unconstrained) assortment.

Using auxiliary variable $ \lambda \in \mathbb{R}_+ $, we can reformulate Problem(\ref{equ:equ4}) as

\begin{equation}\label{equ:equ7}
\begin{aligned}
z^\ast	& =\max_{S,\lambda} \left\{ \lambda :\ \pi(S) \sum_{i \in S} \rho_i(S) r_i \geq \lambda, \ \lambda \in \mathbb{R}_+, \ S \in \varOmega \right\} \\
& =\max_{S,\lambda} \left\{ \lambda :\ \sum_{i \in S} \nu_i (r_i \pi(S)-\lambda) \geq \lambda, \ \lambda \in \mathbb{R}_+, \ S \in \varOmega \right\},	
\end{aligned}
\end{equation}
where the second equality follows from replacing $ \rho_i(S) $ from (\ref{equ:equ00}) and then simplifying the inequality. For any given value of $ \lambda \in \mathbb{R}_+ $, the optimal assortment is given by

\begin{equation}\label{equ:equ8}
\begin{aligned}
S^{\prime}(\lambda) & =\argmax_{S \in \varOmega} \left\{\sum_{i \in S} \nu_i \left( {r_i \pi(S)-\lambda} \right)\right\} ,	
\end{aligned}
\end{equation}

which allows us to reformulate (\ref{equ:equ7}) as

\begin{equation}\label{equ:equ9}
\begin{aligned}
z^\ast & =\max _{\lambda \in \mathbb{R}_+} \left\{\pi(S^{\prime}(\lambda)) \sum_{i \in S^{\prime}(\lambda)} \rho_i(S^{\prime}(\lambda))r_i\right\}.	
\end{aligned}
\end{equation}

We prove that in order to solve unconstrained version of Problem (\ref{equ:equ4}) with $ \bar u = 1 $ it suffices to compute $ S^{\prime}(\lambda) $ for $ O(n) $ values of $ \lambda \in \mathbb{R}_+ $. Furthermore, we show that it is possible to identify these $ O(n) $ values of $ \lambda $, denoted by set $ \Lambda $, in a tractable fashion.

Before describing our solution algorithm, we provide a geometric interpretation of Problem (\ref{equ:equ8}). To do so, we define the linear function $ f_i(\lambda, S)=\nu_i\left( {r_i \pi(S)-\lambda} \right) $ to denote the benefit of product $ i \in S $ for given values of $ \lambda $ and $ S $. Using this definition, we can write the objective function in  (\ref{equ:equ8}) as $ \sum_{i \in S} f_i(\lambda, S) $. Further, when $ \bar u=1 $, we can reformulate $ \pi(S) $ as

\begin{equation}\label{equ:equ14}
\begin{aligned}
\pi(S)=\min_{i \in \bar S} \{1-p^0_i \}.
\end{aligned}
\end{equation}

Because $ f_i(0, S)=\nu_i r_i \pi(S) \geq 0$, for all $  i \in S $ and $ \pi(S) \leq \pi(I \setminus \{0\} )=1 $, for all $ S \in \varOmega $ according to Lemma \ref{l:lambda00}, we have that $ S^{\prime}(\lambda)=I \setminus \{0\} $ when $ \lambda = 0 $. By increasing the value of $ \lambda $, $ f_i(\lambda, I \setminus \{0\}) $ decreases for all $  i \in I \setminus \{0\}$. However, if for a given value of $ \lambda $ we have that $ \sum_ {i \in I} f_i(\lambda, I \setminus \{0\}) > \sum_ {i \in I\setminus \{j, 0\}} f_i(\lambda, I\setminus \{j, 0\}) $ for all $ j \in I \setminus \{0\} $, then $ S^{\prime} (\lambda) = I \setminus \{0\} $. This means that for such value of $\lambda$ it is not optimal to remove any product from the assortment. This result also applies for any other assortment, suggesting a procedure to determine whether removing a product from the assortment is optimal. The following lemma shows that there are at most $ n + 1 $ distinct values from a set $ \Lambda = \{\lambda_0,...,\lambda_{n}\} $ such that $ \lambda > \lambda_k $ implies that it is optimal to remove at least $ k $ products from assortment $ S^{\prime}(0) $.

\begin{lemma} \label{l:lambda}
	Define $ \lambda_0=0 $ and let $ \lambda_k $ be the smallest value of $ \lambda $ that satisfies
	\begin{equation}\label{equ:equ15}
	\begin{aligned}
	\sum_ {i \in S^{\prime}(\lambda_{k-1})} f_i(\lambda_k,  S^{\prime}(\lambda_{k-1})) \leq \sum_ {i \in  S^{\prime}(\lambda_{k-1}) \setminus \{j\}} f_i(\lambda_k,  S^{\prime}(\lambda_{k-1}) \setminus \{j\})
	\end{aligned}
	\end{equation}
	for $ k \in \{1,...,n\} $ and some $ j \in  S^{\prime}(\lambda_{k-1}) $. Then, for any arbitrary assortment $ S \subseteq S^{\prime}(\lambda_{k-1}) $ that includes product $ j $ and for all $ \lambda \geq \lambda_k $, we have
	\begin{equation}\label{equ:equ16}
	\begin{aligned}
	\sum_ {i \in S} f_i(\lambda, S) \leq \sum_ {i \in S\setminus \{j\}} f_i(\lambda, S\setminus \{j\}).
	\end{aligned}
	\end{equation}
	That is, assortment $ S^{\prime}(\lambda) $ does not include product $ j $ for any $ \lambda \geq \lambda_k $.
\end{lemma}

\proof{Proof.}  
First, we show that the set $ \Lambda $ always exists. For any arbitrary value of $ \lambda $ and $ j \in I \setminus \{0\} $ we have that the change in revenue after removing $j$ from assortment $ S^{\prime}(\lambda) $ is given by
\begin{align*}
& \sum_ {i \in S^{\prime}(\lambda)} f_i(\lambda,  S^{\prime}(\lambda)) - \sum_ {i \in  S^{\prime}(\lambda) \setminus \{j\}} f_i(\lambda,  S^{\prime}(\lambda) \setminus \{j\}) = \\
& \sum_{i \in S^{\prime}(\lambda) \setminus \{j\}} \nu_i r_i [ { \pi(S^{\prime}(\lambda))-\pi(S^{\prime}(\lambda) \setminus \{j\})} ]  +  \nu_j r_j \pi(S^{\prime}(\lambda))- \lambda .
\end{align*}

From Lemma \ref{l:lambda00} we have that ${ \pi(S^{\prime}(\lambda))-\pi(S^{\prime}(\lambda) \setminus \{j\})} \geq 0$. Moreover, ${ \pi(S^{\prime}(\lambda))-\pi(S^{\prime}(\lambda) \setminus \{j\})} \leq 1$ and all $\nu$- and $r$-parameters are nonnegative, meaning that there is a value of $ \lambda $, denoted by $ \lambda^\prime_j $, such that the previous expression is negative. This implies that
\begin{align*}
& \sum_ {i \in S^{\prime}(\lambda)} f_i(\lambda,  S^{\prime}(\lambda)) \leq \sum_ {i \in  S^{\prime}(\lambda) \setminus \{j\}} f_i(\lambda,  S^{\prime}(\lambda) \setminus \{j\}), \quad \forall \lambda \geq \lambda^\prime_j .
\end{align*}

Starting with $\lambda_0 = 0$ and applying this argument recursively for $j \in I \setminus \{0\}$, we can build the set $ \Lambda $. Now, define the function
\begin{equation}\label{equ:equ19}
\begin{aligned}
g(S,j)=\frac{\nu_j r_j \pi(S \setminus \{j\})+\sum_ {i \in S} \nu_i r_i (\pi(S)-\pi(S \setminus \{j\}))}{\nu_j}.
\end{aligned}
\end{equation}

Using the definition of $ f_i(\lambda, S) $ and $ f_i(\lambda, S \setminus \{j\}) $, it follows that $ g(S,j) \leq \lambda $ is equivalent to Inequality $ (\ref{equ:equ16}) $. We show that Inequality $ (\ref{equ:equ16}) $ is also valid for every $ \lambda \geq \lambda_k $ and $ S \subseteq S^{\prime}(\lambda_{k-1}) $ that includes $ j $. Let $  \bar S^{\prime}(\lambda_{k-1}) = I \setminus (S^{\prime}(\lambda_{k-1}) \cup \{0\}) $ and suppose that $ \pi(S)=1-p_l^0 $ for some $ \ell \in S $ and that $ \pi(S^{\prime}(\lambda_{k-1}))=1-p_t^0 $ for some $ t \in S^{\prime}(\lambda_{k-1}) $. Since $ \bar S^{\prime}(\lambda_{k-1}) \subseteq \bar S  $, then according to Lemma \ref{l:lambda00} we have $ \pi(S) \leq \pi(S^{\prime}(\lambda_{k-1}))  $ and $ p_t^0 \leq p_\ell^0 $. There are three cases to consider: $ p_j^0 \leq p_t^0 \leq p_\ell^0 $, $  p_t^0 \leq p_j^0 \leq p_\ell^0 $, and $  p_t^0 \leq p_\ell^0 \leq p_j^0 $.

If $ p_j^0 \leq p_t^0 \leq p_\ell^0 $, we have $ \pi(S \setminus \{j\})=\pi(S)=1-p_l^0$ and $ \pi(S^{\prime}(\lambda_{k-1}) \setminus \{j\})=\pi(S^{\prime}(\lambda_{k-1}))=1-p_t^0$. Therefore, it follows that 
\begin{align*}
g(S,j) &  = \frac{\nu_j r_j (1-p_\ell^0)}{\nu_j}
 \leq \frac{\nu_j r_j (1-p_t^0)}{\nu_j}
 = g(S^{\prime}(\lambda_{k-1}),j)
 \leq \lambda.
\end{align*}

Now, suppose that $  p_t^0 \leq p_j^0 \leq p_\ell^0 $. In this case, we have that $ \pi(S \setminus j)=\pi(S)=1-p_\ell^0 $ and $ \pi(S^{\prime}(\lambda_{k-1}) \setminus \{j\})=1-p_j^0 $, which implies
\begin{align*}
g(S,j) &  = \frac{\nu_j r_j (1-p_\ell^0)}{\nu_j} \leq \frac{\nu_j r_j (1-p_j^0)}{\nu_j} \leq \frac{\nu_j r_j (1-p_j^0)+\sum_ {i \in S^{\prime}(\lambda_{k-1})} \nu_i r_i (p_j^0-p_t^0)}{\nu_j} = g(S^{\prime}(\lambda_{k-1}),j)
\leq \lambda.
\end{align*} 

For $  p_t^0 \leq p_\ell^0 \leq p_j^0 $, we have $ \pi(S \setminus \{j\})=1-p_j^0 $ and $ \pi(S^{\prime}(\lambda_{k-1}) \setminus j)=1-p_j^0 $. Then, we have
\begin{align*}
g(S,j) &  = \frac{\nu_j r_j (1-p_j^0)}{\nu_j}
 = g(S^{\prime}(\lambda_{k-1}),j)
 \leq \lambda.
\end{align*} 

In any case we obtain that $g(S,j)
 \leq \lambda$, which implies that $\sum_ {i \in S} f_i(\lambda, S) \leq \sum_ {i \in S\setminus \{j\}} f_i(\lambda, S\setminus \{j\})$, completing the proof. $ \square  $

Proposition \ref{p:LineMethod} follows immediately from Lemma \ref{l:lambda} and the definition of $ \Lambda$.

\begin{proposition} \label{p:LineMethod} $ S^{\prime}(\lambda) \subseteq S^{\prime}(\lambda_{k-1}) \setminus \{j\} $ for any $ \lambda \geq \lambda_k $, all $ k \in \{1,...,n\} $, and some $ j \in S^{\prime}(\lambda_{k-1}) $. 
\end{proposition} 

\proof{Proof.} From Lemma \ref{l:lambda}, we know that $S^{\prime}(\lambda) \subseteq S^{\prime}(\lambda_k) \setminus \{j\}$ and that $S^{\prime}(\lambda_k) \setminus \{j\} \subseteq S^{\prime}(\lambda_{k-1}) \setminus \{j\}$ for $ \lambda \geq \lambda_k $ and some $j \in S^{\prime}(\lambda_{k-1})$. This implies that $ S^{\prime}(\lambda) \subseteq S^{\prime}(\lambda_{k-1}) \setminus \{j\} $, proving the result. $ \square  $

This proposition implies that the optimal solution to the unconstrained \APP{} with $ \bar u = 1 $ can be found in polynomial time by generating at most $ n + 1 $ candidate assortments, which can be evaluated based on their objective function value in Problem (\ref{equ:equ9}). We emphasize that set $\Lambda$ has at most $n+1$ elements as some $\lambda$-values may be repeated. Our polynomial-time solution method to solve Problem $ (\ref{equ:equ9}) $ is formally presented in Algorithm \ref{alg:1}. In this algorithm, we find the candidate assortments $ S^\prime(\lambda) $ for all $ \lambda \in \Lambda$ and select the one with the maximum objective value in Problem $ (\ref{equ:equ9}) $. In Line $1$, we initialize $\lambda_0 = 0$, $ S^{\prime}(\lambda_0) = I \setminus \{0\} $, $ S^\ast = S^\prime(\lambda_0) $ and its corresponding objective value $ z^\ast $. Using the results from Lemma \ref{l:lambda} and Proposition \ref{p:LineMethod}, the loop in Lines $2-14$ finds the elements of $\Lambda$ and their corresponding assortment $ S^\prime(\lambda_k) $ for all $ k \in \{1,...,n\} $. In line $3$, we initialize $ \lambda_k = \infty$ and in Line $4$, we find the value of $ \pi({S^{\prime}}(\lambda_{k-1})) $ which is used to evaluate the $g$-functions later on. In Lines $ 5-11 $, we test whether removing Product $j$ from assortment $ S^\prime(\lambda_{k-1}) $ can be used to update $ \lambda_{k} $ and $ S^\prime(\lambda_{k}) $. In Line $6$, we calculate $ g(S^{\prime}(\lambda_{k-1}),j) $, and in Lines $7-8$, we use the condition $ g(S^{\prime}(\lambda_{k-1}),j) \leq \lambda_k $ as it is equivalent to Inequality (\ref{equ:equ15}). These lines determine the product to remove and its corresponding $\lambda$-value, where $ j^\prime $ is such that $ S^{\prime}(\lambda_k)=S^{\prime}(\lambda_{k-1}) \setminus \{j^{\prime}\} $. Note that if $\lambda_k \leq \lambda_{k-1}$, it is also optimal to remove product $j^\prime$ from assortment $S^\prime(\lambda_{k-1})$. Therefore, Lines $9-10$ verify this condition and set $\lambda_k = \lambda_{k-1}$ to be able to update assortment $S^\prime(\lambda_{k-1})$ in Line $11$. In Lines $12-14$, we update $ S^\ast $ and $z^\ast$ if the assortment $ S^\prime(\lambda_k) $ in the $ k $th iteration is more profitable than the best observed thus far. The complexity of Algorithm \ref{alg:1} is $O(n^2)$ due to the for-loops starting in Lines $ 2 $ and $ 5 $. The correctness of Algorithm \ref{alg:1} follows from Lemma \ref{l:lambda} and Proposition \ref{p:LineMethod}. 

\begin{algorithm} 
	\footnotesize{
		\caption{: Solution algorithm for unconstrained \APP{} with $ \bar u = 1 $ \label{alg:1}}
		\begin{algorithmic}[1]
			\State Initialize $ \lambda_0 =0 $, $ S^{'}(\lambda_0)=I \setminus \{0\} $, $ S^\ast=S^{'}(\lambda_0) $, and $ z^\ast=\sum_{i \in S^{'}(\lambda_0)} \rho_i(S^{'}(\lambda_0))r_i $
			\ForAll{$k \in  \{1,...,n\}$}
			    \State Initialize $ \lambda_k =\infty $
			    \State Set $ \pi({S^{'}}(\lambda_{k-1}))=\min_{i \in \bar S^{'}(\lambda_{k-1})} \{1-p^0_i \} $
		    	\ForAll {$j \in {S^{'}}(\lambda_{k-1}) $}
			        \State Set $ \pi({S^{'}}(\lambda_{k-1}) \setminus \{j\})=\min_{i \in \bar S^{'}(\lambda_{k-1}) \cup \{j\} } \{1-p^0_i \} $ and calculate $ g(S^{'}(\lambda_{k-1}),j) $ using Equation (\ref{equ:equ19})
		    	    \If {$ g(S^{'}(\lambda_{k-1}),j) \leq \lambda_k $} 
		        	    \State Set $ j^{'}=j $ and $ \lambda_k=g(S^{'}(\lambda_{k-1}),j) $
		    	    \EndIf
		        \EndFor
		    	\If {$ \lambda_k \leq \lambda_{k-1}$} 
		    	    \State Set $\lambda_k = \lambda_{k-1}$
		    	\EndIf
		    	\State Set $ S^{'}(\lambda_k)=S^{'}(\lambda_{k-1}) \setminus \{j^{'}\} $  
		    	\If {$ z^\ast \leq \pi(S^{'}(\lambda_k)) \sum_{i \in S^{'}(\lambda_k)} \rho_i(S^{'}(\lambda_k))r_i $ }   
		        	\State Set $ z^\ast = \pi(S^{'}(\lambda_k)) \sum_{i \in S^{'}(\lambda_k)} \rho_i(S^{'}(\lambda_k))r_i $ 
		    	    \State $ S^\ast= S^{'}(\lambda_k)$     
		    	\EndIf
			\EndFor
			\State \textbf{return} $ S^\ast$
	\end{algorithmic}}
\end{algorithm}
\normalsize        

The following example illustrates the operation of Algorithm \ref{alg:1} in a problem with $3$ candidate products.

\begin{figure}[htbp]
	\begin{center}
		\centering
		\includegraphics[scale=0.48]{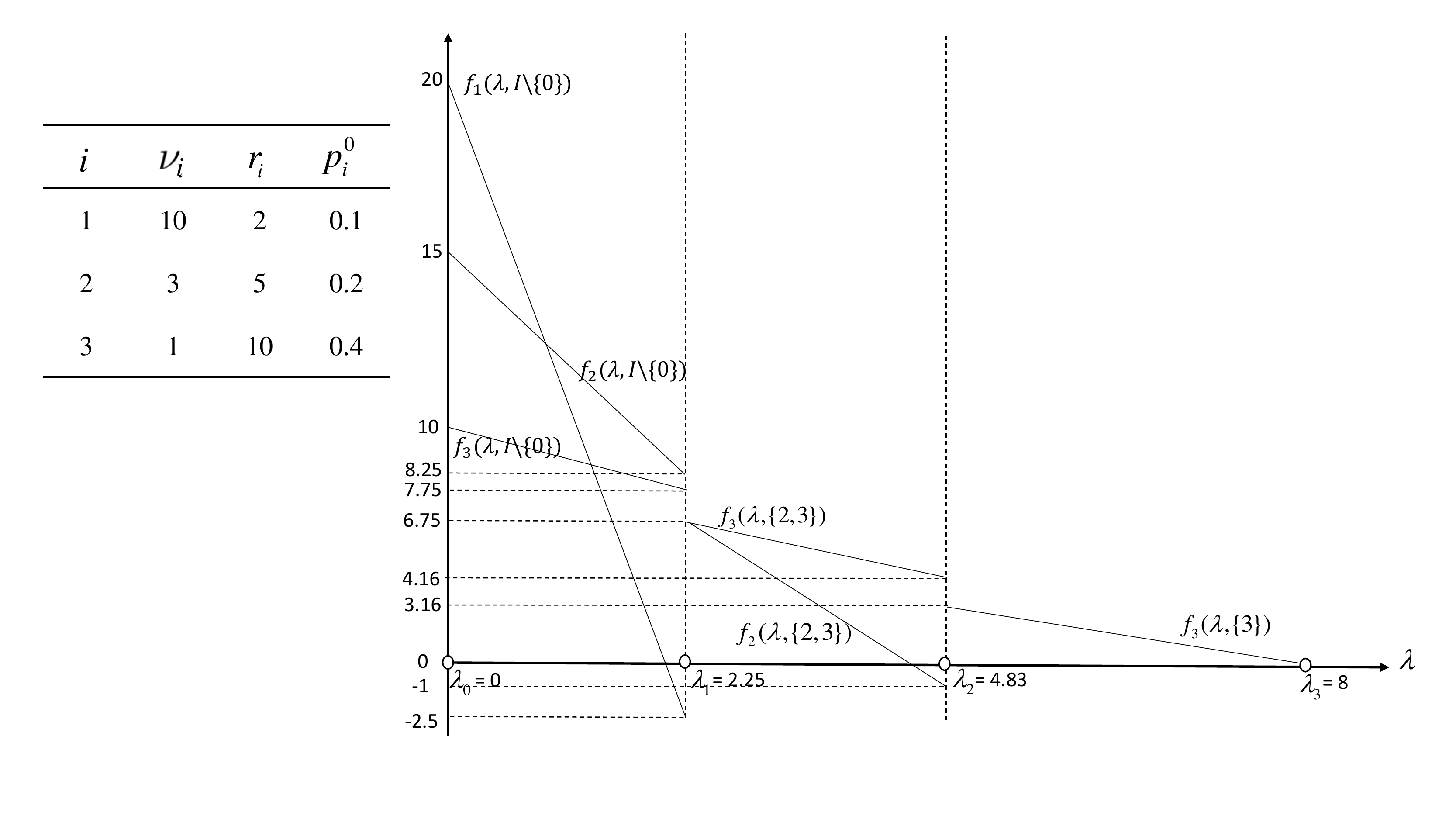}
		\caption[]{An illustration of our results in Example \ref{E:LineMethod} with three products}
		\label{f:Example2}
	\end{center}
\end{figure}

\begin{example} \label{E:LineMethod} Suppose that there are $3$ products available with parameters and functions $ f_i(S^{\prime}(\lambda), \lambda) $ given by Figure $ \ref{f:Example2} $. For $ \lambda=0 $, we have that $ S^{\prime}(\lambda)=I \setminus \{0\} $ because $ f_i(I \setminus \{0\}, \lambda) \geq 0 $ for every product. However, as $\lambda$ increases, the profits of products become negative one by one, starting with Product $1$ at $\lambda=\lambda_1$. For any $ \lambda \geq \lambda_1 $, we have $\sum_{i \in I \setminus \{0\}} f_i(I \setminus \{0\}, \lambda) \leq \sum_{i \in \{2,3\}} f_i(\{2,3\}, \lambda)$, indicating that Product $1$ must be removed from the assortment at the expense of decreasing $ \pi(\cdot) $. Accordingly, we obtain $ S^{\prime}(\lambda)= \{2,3\} $ for $ \lambda \in [\lambda_1, \lambda_2) $. Using the same argument, Product $ 2 $ is not in the optimal assortment $S^\prime(\lambda) = \{3\}$ for $ \lambda \in [\lambda_2, \lambda_3)$. Finally, for $ \lambda \geq \lambda_3 $, all the three products have negative profit and $ S^{\prime}(\lambda)=\emptyset $. As a result, the candidate assortments corresponding to $\lambda$-segments $[\lambda_{k-1},\lambda_k)$ for $k =1,...,4$, where $\lambda_4 = \infty$, are $ I \setminus \{0\}$ ($z=3$), $ \{2,3\} $ ($z=4.5$), $ \{3\} $ ($z=4$), and $ \emptyset $ ($z=0$), respectively, where $ S=\{2,3\} $ is the optimal solution to Problem $ (\ref{equ:equ9}) $.
\end{example}

\subsection{Algorithms for constrained \APP{}} \label{subsec:ConSolution} 
In this section, we generalize Problem (\ref{equ:equ4}) by assuming that the maximum number of products allowed in the assortment is $ \bar c $, i.e., $ \varOmega= \left\{S: S \subseteq I, |S| \leq \bar c \right\} $. Before presenting our solution approaches, the following proposition provides a reformulation of Problem $ (\ref{equ:equ3}) $.

\begin{proposition} \label{p:LinearProgram}
	Problem $ (\ref{equ:equ3}) $ can be transformed into a linear program to find $ \pi(S) $ for any assortment $ S \subseteq I \setminus \{0\} $.
\end{proposition}

\proof{Proof.} Because of the objective function of (\ref{equ:equ3}) and the nonnegativity of $\eta$- and $p^0$-parameters, we can replace $\sum_{i \in \bar S}y_{ik}= 1$ with $\sum_{i \in \bar S}y_{ik}\leqslant 1$ for each $k \in \{1,...,{\bar u}\}$. The resulting constraint coefficient matrix is totally unimodular, which allows us to reformulate $ Y(S) $ as
\begin{align*}
Y^\prime(S)= \left\{ \mathbf{y} \in [0,1] ^ {{\lvert \bar S \rvert} \times {\bar u}} :\ \sum_{i \in \bar S}y_{ik}\leqslant 1, \ \forall k \in \{1,...,{\bar u}\}, \  \sum_{k=1}^{\bar u}y_{ik}\leqslant 1, \ \forall i \in \bar S \right\},
\end{align*}
whose extreme points are integer-valued. Because the natural logarithm is a monotonically increasing function and $ \sum_{k=1}^{\bar u}y_{ik}\leqslant 1$, $ \forall i \in \bar S $, we can transform Problem $ (\ref{equ:equ3}) $ into the linear program
\begin{equation}\label{equ:equ6}
\begin{aligned}
\theta(S) &= \min_{\mathbf{y} \in Y^\prime(S)} \ln \left(\prod_{k=1}^{{\bar u}}\left(1-\sum_{i \in \bar S} {\eta_{ik} p_i^0 y_{ik}}\right)\right) \\
&= \min_{\mathbf{y} \in Y^\prime(S)} \sum_{k=1}^{{\bar u}} \ln \left(1-\sum_{i \in {\bar S}} {y_{ik} \eta_{ik} p^0_{i}}\right) \\
& =\min_{\mathbf{y} \in Y^\prime(S)} \sum_{k=1}^{{\bar u}} \sum_{i \in {\bar S}} y_{ik}  \ln (1-{\eta_{ik} p^0_{i}}),
\end{aligned}
\end{equation}
where the first and second equalities follow from the properties of the natural logarithm function and the third equality holds because $ \sum_{i \in {\bar S}} y_{ik} $ is either equal to $0$ or $1$ for any $k \in \{1,...,{\bar u}\}$. Therefore, $ \theta(S) $ can be found by solving the linear program $ (\ref{equ:equ6}) $ from which we can calculate  $ \pi(S)=e^{\theta(S)} $, proving the proposition. $ \square $ 

Using the results from Proposition \ref{p:LinearProgram}, we transform the bi-level Problem $ (\ref{equ:equ4}) $ into a single-level problem. For any assortment $ S $, the dual problem of the linear program $ (\ref{equ:equ6}) $ is given by
\begin{equation}\label{equ:equ21}
\begin{aligned}
\max_{\boldsymbol{\alpha}, \boldsymbol{\beta} \in A(S)} \left\{-\sum_{k=1}^{\bar u} \alpha_k-\sum_{i \in \bar S}\beta_i \right\},
\end{aligned}
\end{equation}
where vectors $ \boldsymbol{\alpha} $ and $ \boldsymbol{\beta} $ are dual decision variables constrained by
\begin{align*}
A(S)= \left\{\boldsymbol{\alpha} \in \mathbb{R}_+^{\bar u}, \boldsymbol{\beta} \in \mathbb{R}_+^{|\bar S|} :\ \alpha_k+\beta_i \geq -\ln(1-\eta_{ik}p_{i0}),\ \forall i \in \bar S,\ k \in \{1,...,\bar u\}\right\}. 
\end{align*}

Using $ \pi(S)=e^{\theta (S)} $ and Proposition \ref{p:LinearProgram}, Problem $ (\ref{equ:equ4}) $ can be reformulated as
\begin{equation}\label{equ:equ23}
\begin{aligned}
z^\ast= & \max_{S \in \varOmega} \left\{ e^{\theta (S)} \sum_{i \in S} \rho_i(S)r_i \right\}.
\end{aligned}
\end{equation}

Combining Problems $ (\ref{equ:equ6}) $ and $ (\ref{equ:equ21}) $, and using the strong duality theorem from linear programming, we obtain that
\begin{equation}\label{equ:equ24}
\begin{aligned}
\theta(S)=\max_{\boldsymbol{\alpha}, \boldsymbol{\beta} \in A(S)} \left\{-\sum_{k=1}^{\bar u} \alpha_k-\sum_{i \in \bar S}\beta_i \right\}.
\end{aligned}
\end{equation}

Because the orientation of both Problems $ (\ref{equ:equ23}) $ and $ (\ref{equ:equ24}) $ is maximization and $ e^{\theta(S)} $ is a monotonically increasing function on $ \theta(S) $, Problem (\ref{equ:equ23}) can be transformed into
\begin{equation}\label{equ:equ30}
\begin{aligned}
z^\ast & = \max_{S \in \varOmega, \boldsymbol{\alpha}, \boldsymbol{\beta} \in A(S)} \left\{ e^{(-\sum_{k=1}^{\bar u} \alpha_k-\sum_{i \in \bar S}\beta_i)} \sum_{i \in S} \rho_i(S)r_i \right\}.
\end{aligned}
\end{equation}

To solve the single-level Problem $ (\ref{equ:equ30}) $, we propose a problem-specific exact cutting-plane algorithm and a greedy solution approach, which are discussed in Sections \ref{subsubsec:Piecewise} and \ref{subsubsec:Greedy}, respectively.

\subsubsection{Exact cutting-plane algorithm}$  $\\
\label{subsubsec:Piecewise} 
We define parameters $ \theta_{min} $ and $ \theta_{max} $ such that $\theta_{min} \leq \min_{S \in \varOmega} \theta(S)$ and $ \theta_{max} \geq \max_{S \in \varOmega} \theta(S)$. Because $\theta(S) \leq 0 $ holds for all $ S \in \varOmega $, we use $ \theta_{max}=0 $. Moreover,  $ \bar{S}=(I \setminus S) \subseteq I $ for any $ S \in \varOmega $ and then we have that $ \theta(\emptyset) \leq \min_{S \in \varOmega} \theta(S)$. Therefore, we use $ \theta_{min}=\theta(\emptyset) $, which can be found by solving Problem $ (\ref{equ:equ6}) $.

Consider the graph of the function $ \phi (\theta(S) ) = e^{\theta(S)}$, for which  $ \phi (\theta(S) ) = \pi (S) $ for any $ S \subseteq I \setminus \{0\} $. Because $ \phi (\theta(S) ) $ is a convex function, the line segment connecting points $ (\theta_{min}, e^{\theta_{min}}) $ and $ (\theta_{max}, e^{\theta_{max}}) $ lies above $ \phi (\theta(S) ) $. Using this observation, we define the function $ h(\theta(S)) $ that provides an upper bound to $ \phi (\theta(S) )$.
\begin{align*}
h(\theta(S))= & \frac{\phi (\theta_{max} )-\phi (\theta_{min} )}{\theta_{max}-\theta_{min}} \left(\theta(S)-\theta_{min}\right)+\phi (\theta_{min} ).
\end{align*} 

The following proposition presents a strategy to construct an upper bound for $ z^{*} $ using $ h(\theta(S)) $ instead of $ e^{\theta(S)} $ in (\ref{equ:equ30}).

\begin{proposition} \label{p:UpperBound} 
	The problem
	\begin{equation}\label{equ:equ22}
	\begin{aligned}
	\bar z= \max_{S \in \varOmega,\ \boldsymbol{\alpha}, \boldsymbol{\beta}  \in A(S)} \left\{h \left(-\sum_{k=1}^{\bar u} \alpha_k-\sum_{i \in \bar S}\beta_i \right) \sum_{l \in S} \rho_l(S)r_l \right\},
	\end{aligned}
	\end{equation}
provides an optimal objective function value $ 	\bar z $ such that $ z^\ast \leq \bar z $.
\end{proposition}
\proof{Proof.}
We have that $ e^{(-\sum_{k=1}^{\bar u} \alpha_k-\sum_{i \in \bar S}\beta_i)} \leq h(-\sum_{k=1}^{\bar u} \alpha_k-\sum_{i \in \bar S}\beta_i) $ for any $S \in\Omega$, and the result follows by replacing $e^{(-\sum_{k=1}^{\bar u} \alpha_k-\sum_{i \in \bar S}\beta_i)}$ with $h(-\sum_{k=1}^{\bar u} \alpha_k-\sum_{i \in \bar S}\beta_i)$ in Problem (\ref{equ:equ30}). $ \square $

\begin{corollary} \label{cor:UpperBound} 
Let $(S^\ast , \boldsymbol{\alpha}^\ast, \boldsymbol{\beta}^\ast)$ be an optimal solution to Problem (\ref{equ:equ22}). Then $ \theta(S^\ast) = -\sum_{k=1}^{\bar u} \alpha^\ast_k-\sum_{i \in \bar S}\beta^\ast_i $.
\end{corollary}

\proof{Proof.}
Decision variables $  \boldsymbol{\alpha}$ and $ \boldsymbol{\beta}$ in Problem (\ref{equ:equ22}) only appear in $ h(\cdot) $ and they are not part of the feasible set $ \Omega $. Therefore, we have
\begin{align*}
\bar z &= \max_{S \in \varOmega} \left\{ \left(\max_{ \boldsymbol{\alpha},  \boldsymbol{\beta} \in A(S)} h \left(-\sum_{k=1}^{\bar u} \alpha_k-\sum_{i \in \bar S}\beta_i \right)\right) \sum_{l \in S} \rho_l(S)r_l \right\} \\
& = \max_{S \in \varOmega} \left\{ \left(\max_{ \boldsymbol{\alpha},  \boldsymbol{\beta} \in A(S)} \left( \frac{\phi (\theta_{max} )-\phi (\theta_{min} )}{\theta_{max}-\theta_{min}} \left(-\sum_{k=1}^{\bar u} \alpha_k-\sum_{i \in \bar S}\beta_i -\theta_{min}\right)+\phi (\theta_{min}) \right) \right) \sum_{l \in S} \rho_l(S)r_l \right\},
\end{align*}
where the second equality results by replacing $ h(\cdot) $ by its definition. Because $ \frac{\phi (\theta_{max} )-\phi (\theta_{min} )}{\theta_{max}-\theta_{min}} \geq 0 $ and $\phi (\theta_{min} )$ is a constant, we obtain the equivalent problem
\begin{align} \label{ubreformulation}
\bar z &= \max_{S \in \varOmega} \left\{ \left( \frac{\phi (\theta_{max} )-\phi (\theta_{min} )}{\theta_{max}-\theta_{min}} \left(\max_{\boldsymbol{\alpha}, \boldsymbol{\beta} \in A(S)} \left( -\sum_{k=1}^{\bar u} \alpha_k-\sum_{i \in \bar S}\beta_i \right) -\theta_{min}\right)+\phi (\theta_{min}) \right) \sum_{l \in S} \rho_l(S)r_l \right\}.
\end{align}

Suppose for a contradiction that the given optimal solution to Problem (\ref{equ:equ22}) is such that $ \theta(S^\ast) > -\sum_{k=1}^{\bar u} \alpha^\ast_k-\sum_{i \in \bar S}\beta^\ast_i $. This means that the optimal value of the inner maximization problem in (\ref{ubreformulation}) is $-\sum_{k=1}^{\bar u} \alpha^\ast_k-\sum_{i \in \bar S}\beta^\ast_i < \theta(S^\ast)$. However, this inner maximization problem is the same as Problem (\ref{equ:equ24}) with $S=S^\ast$, which implies that $\boldsymbol{\alpha}^\ast$ and $ \boldsymbol{\beta}^\ast $ cannot be optimal because the optimal solution to Problem (\ref{equ:equ24}) (with value $\theta(S^\ast)$) is feasible to Problem (\ref{equ:equ22}). Because $ \theta(S^\ast) < -\sum_{k=1}^{\bar u} \alpha^\ast_k-\sum_{i \in \bar S}\beta^\ast_i $ contradicts the optimality of $\theta(S^\ast)$, we conclude that $ \theta(S^\ast) = -\sum_{k=1}^{\bar u} \alpha^\ast_k-\sum_{i \in \bar S}\beta^\ast_i $, completing the proof. $ \square $

Let $ \mathbf{x}=\{x_1,...,x_n\} \in \{0,1\}^n $ be a vector of binary decision variables, where $ x_i $ equals one if product $ i $ is included in the assortment and equals zero, otherwise, such that $S = \{i \in I \setminus \{0\}: x_i=1\}$. Accordingly, we rewrite the feasible region $A(S)$ as
\begin{align*}
A(\mathbf{x})= \left\{\boldsymbol{\alpha} \in \mathbb{R}_+^{\bar u}, \boldsymbol{\beta} \in \mathbb{R}_+^{n} :\ \alpha_k+(1-x_i)\beta_i \geq -\ln(1-\eta_{ik}p_{i0})(1-x_i),\ \forall i=1,\dots,n ,\ k \in \{1,...,\bar u\}\right\} \end{align*}
and define the assortment feasible region $X=\{\mathbf{x} \in \{0,1\}^n: \sum_{i=1}^n x_i \leq \bar c\}$, which allow us to rewrite Problem $ (\ref{equ:equ22}) $ as
\begin{align*}
\bar z=\max_{\mathbf{x} \in X,\ \boldsymbol{\alpha}, \boldsymbol{\beta} \in A(\mathbf{x})} \left\{ h \left(-\sum_{k=1}^{\bar u} \alpha_k-\sum_{i =1}^n(1-x_i)\beta_i \right) \sum_{i =1}^n \frac{\nu_i r_i x_i}{1+\sum_{j =1}^n \nu_j x_j} \right\} .
\end{align*}

Using linear fractional programming techniques, we can further reformulate the problem by introducing the auxiliary variable $ \delta \in \mathbb{R} $.
\begin{equation}\label{equ:equ31}
\begin{aligned}
\bar z &=\max \left\{ \delta:\  h \left(-\sum_{k=1}^{\bar u} \alpha_k-\sum_{i=1}^n(1-x_i)\beta_i \right) \sum_{i =1}^n \frac{\nu_i r_i x_i}{1+\sum_{j =1}^n \nu_j x_j} \geq \delta, \mathbf{x} \in X,\ \boldsymbol{\alpha}, \boldsymbol{\beta} \in A(\mathbf{x}),\ \delta \in \mathbb{R} \right\} \\
& = \max \left\{ \delta:\ \sum_{j =1}^n {\nu_j x_j \left(r_j h \left(-\sum_{k=1}^{\bar u} \alpha_k-\sum_{i =1}^n(1-x_i)\beta_i \right)-\delta \right)} \geq \delta,\ \mathbf{x} \in X,\ \boldsymbol{\alpha}, \boldsymbol{\beta} \in A(\mathbf{x}),\ \delta \in \mathbb{R} \right\},
\end{aligned}
\end{equation}
where the second equality follows from multiplying both sides of the inequality by the denominator of the fraction and then reorganizing the expression. This problem can be linearized using conventional techniques \citep{gupte2013solving}. As a result, the optimal assortment to Problem (\ref{equ:equ31}) can be found by solving a mixed-integer program (\MIP{}), providing an upper bound on the optimal value to Problem (\ref{equ:equ23}). 

\begin{figure}[htbp]
	\begin{center}
		\centering
		\includegraphics[scale=0.5]{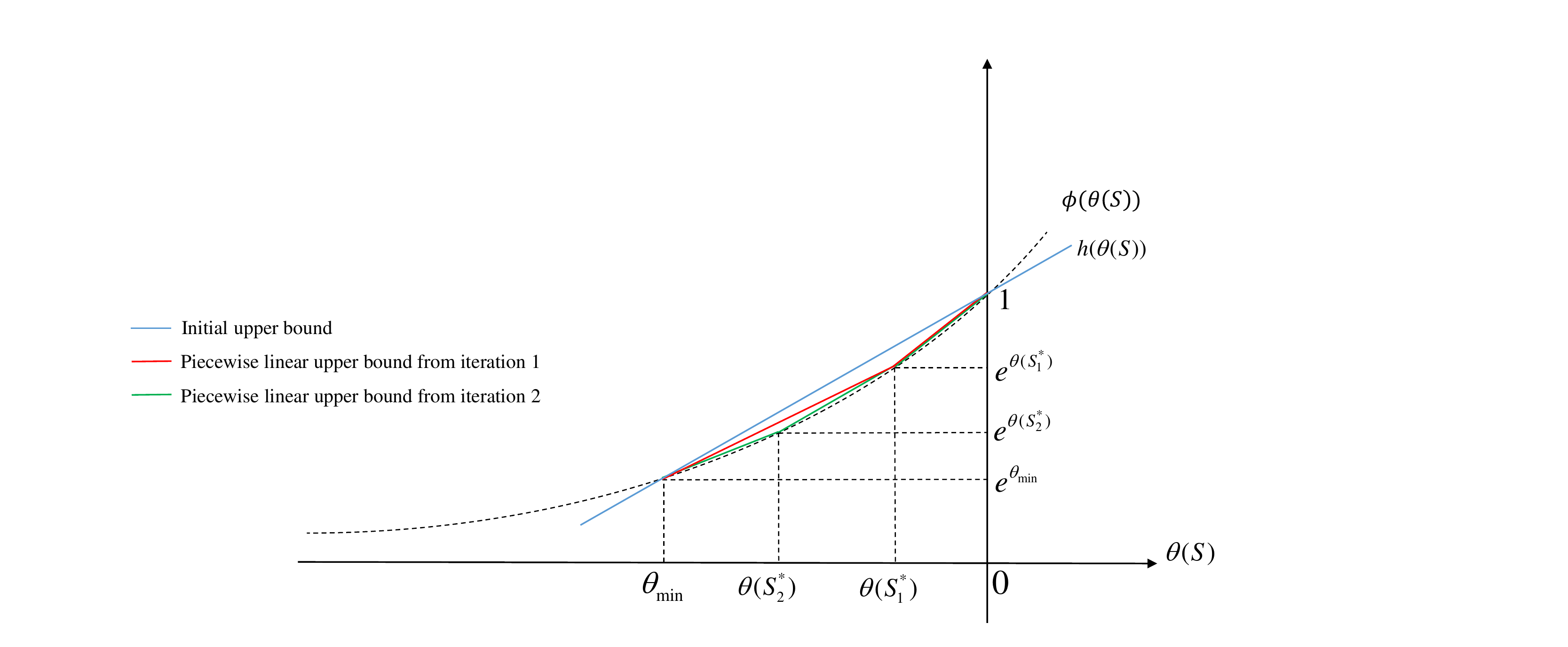}
		\caption[]{Piecewise linear updated function $ h(\theta(S)) $ based on assortments $ S^\ast_1 $ and $ S^\ast_2 $}
		\label{f:3}
	\end{center}
\end{figure}

We iteratively tighten the upper-bound Problem (\ref{equ:equ31}) by transforming the linear function $ h(\cdot) $ into a piecewise-linear function based on the assortments found at each iteration. To illustrate this approach, suppose that Problem (\ref{equ:equ31}) is solved, producing an optimal assortment $ S^\ast_1 $, with $ h(\theta (S^\ast_1))$ calculated by $ h \left(-\sum_{k=1}^{\bar u} \alpha^\ast_k-\sum_{i \in I}(1-x^\ast_i)\beta^\ast_i \right)$ given the result in Corollary \ref{cor:UpperBound}. If $ \phi (\theta(S^\ast_1) ) < h(\theta(S^\ast_1)) $, then we refine the piecewise-linear approximation by adding two lines with common endpoint $ (\theta(S^\ast_1), \phi (\theta(S^\ast_1) )) $, which produces the exact value of $ \theta(S) $ when assortment $ S^\ast_1 $ is selected. Because $ \phi (\theta(S^\ast_1) ) < h(\theta(S^\ast_1)) $, these new line segments are guaranteed to produce a tighter upper bound on $ \phi (\theta(S^\ast_1) ) $. Problem (\ref{equ:equ31}) is then solved again to obtain a new assortment and the tightening procedure is applied again, if needed. Figure \ref{f:3} illustrates the operation of our approach for two iterations, assuming that assortments $ S^\ast_1 $ and $ S^\ast_2 $ are obtained in the first and the second iteration, respectively. We refer to this approach as \emph{piecewise-linear upper bound} (\PLUB{}).

Formally, let $ T $ be the number of times that Problem (\ref{equ:equ31}) has been solved using \PLUB{}. Suppose that $ \{S^\ast_1,...,S^\ast_T\} $ is the set of assortments obtained thus far such that $ \theta_{min} \leq \theta(S^\ast_1) \leq ... \leq \theta(S^\ast_T) \leq 0 $. Also, suppose that $ \phi (\theta(S^\ast_T) ) < h(\theta(S^\ast_T))$ so that the upper bound approximation can be improved. Define auxiliary variables $ \mathbf{w}=[w_1,...,w_{T+2}] \in [0,1]^{T+2} $, $ \mathbf{q}=[q_1,...,q_{T+1}] \in \{0,1\}^{T+1} $, and $\pi^\prime \in \mathbb{R} $. Thus, the piecewise-linear approximation used at iteration $ T+1 $ is defined by the following set of constraints, which are added to Problem $ (\ref{equ:equ31}) $.
\begin{subequations}
\begin{align}
 & w_1 \leq  q_1  \label{equ:piecewise1}\\
 & w_t \leq  q_{t-1}+q_t, \quad \forall t \in \{2,...,T+1\}  \label{equ:piecewise2}\\
 & w_{T+2} \leq  q_{T+1} \label{equ:piecewise3}\\
 & \sum_{t=1}^{T+1} q_t = 1 \label{equ:piecewise4}\\
 & \sum_{t=1}^{T+2} w_t = 1 \label{equ:piecewise5}\\
 & \sum_{t=1}^{T+2} \theta(S^{\ast}_t) w_t =-\sum_{k=1}^{\bar u} \alpha_k-\sum_{i=1}^n(1-x_i)\beta_i \label{equ:piecewise6}\\
 & \pi^\prime \leq \sum_{t=1}^{T+2} e^{\theta(S^{\ast}_t)} w_t \label{equ:piecewise7}\\
 & \mathbf{w} \in [0,1]^{T+2}, \mathbf{q} \in \{0,1\}^{T+1}, \pi^\prime \in \mathbb{R}. \label{equ:piecewise8}
\end{align}
\end{subequations}

Constraints (\ref{equ:piecewise1})--(\ref{equ:piecewise8}) describe a conventional \textit{convex combination} formulation to model piece-wise linear function, where continuous $w$-variables are the weights used to construct any value of $\theta(S)$ (i.e., $x$-axis domain) as a convex combination of the endpoints of consecutive line segments and $q$-variables control which segments are used \citep{vielma2010mixed}. In particular, Constraints (\ref{equ:piecewise1})-(\ref{equ:piecewise3}) are for weight activation and Constraint  (\ref{equ:piecewise4}) is for selecting only one segment. Constraints (\ref{equ:piecewise5}) and (\ref{equ:piecewise6}) generate the convex combination of consecutive endpoints. Constraint (\ref{equ:piecewise7}) allows us to replace $ h \left(-\sum_{k=1}^{\bar u} \alpha_k-\sum_{i \in I}(1-x_i)\beta_i \right) $ in Problem $ (\ref{equ:equ31}) $ with $ \pi^\prime $ at iteration $ T $. Constraints (\ref{equ:piecewise8}) enforce the nature of the decision variables. We define the feasible set produced by Constraints (\ref{equ:piecewise1})--(\ref{equ:piecewise8}) as $ Q_T $. Therefore, we can rewrite Problem $  (\ref{equ:equ31}) $ at iteration $ T $ as  
\begin{equation}\label{equ:equ32}
\begin{aligned}
\bar z_T & = \max \left\{ \delta:\ \sum_{j =1}^n \left({\nu_j x_j \left(r_j \pi^\prime -\delta \right)} \right)\geq \delta, \mathbf{x} \in X,\ \boldsymbol{\alpha}, \boldsymbol{\beta} \in A(\mathbf{x}),\ \delta \in \mathbb{R},\ (\mathbf{q}, \mathbf{w}, \pi^\prime) \in Q_T \right\},
\end{aligned}
\end{equation}
where the nonlinear terms $ \pi^\prime x_j $ and $ x_j \delta $ can be linearized using conventional techniques. 

Algorithm \ref{alg:2} describes the steps performed by our \PLUB{} approach. We initialize $ T $, $ \theta_{min} $, and $ \theta_{max} $ in Line $ 1 $. In Line $ 2 $, we calculate an initial assortment by solving Problem (\ref{equ:equ31}) and calculate its corresponding $ \theta(S^\ast_T)$ and $ \pi^\prime $ values in Line $ 3 $. The while-loop in Lines $4$--$8$ verifies the condition $\pi ^ \prime > \phi (\theta(S^\ast_T))$, which indicates, if satisfied, that the upper bound approximation still overestimates $\pi(S)$. If this condition is not satisfied, we update $ T $ and  $ Q_T $ and solve Problem (\ref{equ:equ32}). By solving Problem (\ref{equ:equ32}) in Line $ 7 $, we always maintain an upper bound for $ \phi(\theta(S)) $ for any assortment $ S \in \Omega $. Moreover, the construction in Line $ 6 $ guarantees that the solution $ (S^\ast_T, \pi^\prime) $ at iteration $ T $ becomes infeasible at iteration $ T+1 $. Algorithm \ref{alg:2} terminates if $ \pi ^ \prime = \phi (\theta(S^\ast_T)) $, which can be achieved if $ S^\ast_T \equiv S^\ast_{T-1} $ or $\theta(S^\ast_T) \equiv \theta(S^\ast_{T-1})$ at any iteration $ T $. Therefore, the worst-case number of iterations in Algorithm \ref{alg:2} is equal to the number of feasible assortments, meaning that the algorithm terminates. The correctness of Algorithm \ref{alg:2} comes from the fact that the solution to Problem (\ref{equ:equ32}) provides an upper bound on $z^\ast$, whereas the calculation of $\phi(\theta(S^\ast_T))$ provides a lower bound on $z^\ast$ given that $S^\ast_T$, $\mathbf{x}^\ast$, $ \boldsymbol{\alpha}^\ast$, and $ \boldsymbol{\beta}^\ast $ are feasible to Problem (\ref{equ:equ30}). Therefore, when $\pi ^ \prime _T = \phi (\theta(S^\ast_T))$ then the upper bound solution to Problem (\ref{equ:equ32}) is feasible to Problem (\ref{equ:equ30}) so it must be optimal.  

\begin{algorithm} 
	\small{
		\caption{: \PLUB{} approach for the constrained \APP{}\label{alg:2}}
		\begin{algorithmic}[1]
			\State Set $ T =0 $, $ \theta_{min} = \theta(\emptyset) $, and $ \theta_{max} = 0 $
			\State Solve upper bound Problem (\ref{equ:equ31}) to obtain an initial optimal assortment $ S^\ast_T $ and $\mathbf{x}^\ast$, $ \boldsymbol{\alpha}^\ast$, and $ \boldsymbol{\beta}^\ast $ values
			\State Set $ \theta(S^\ast_T) = -\sum_{k=1}^{\bar u} \alpha^\ast_k-\sum_{i =1}^n (1-x^\ast_i)\beta^\ast_i  $ and $ \pi^\prime = h(\theta(S^\ast_T)) $
			\While {$\pi ^ \prime  - \phi (\theta(S^\ast_T)) > 0$}
			\State Set $ T = T + 1 $ 
			\State Construct set of constraints $ Q_T $
			\State Solve upper bound Problem (\ref{equ:equ32}) to obtain an optimal assortment $ S^\ast_{T} $ and $\mathbf{x}^\ast$, $ \boldsymbol{\alpha}^\ast$, and $ \boldsymbol{\beta}^\ast $ values
			\State Set $ \theta(S^\ast_{T}) = -\sum_{k=1}^{\bar u} \alpha^\ast_k-\sum_{i =1}^n (1-x^\ast_i)\beta^\ast_i  $
			\EndWhile
			\State \textbf{return} $ S_T^\ast$
	\end{algorithmic}}
\end{algorithm}
\normalsize  

We further tighten the formulation of Problem (\ref{equ:equ32}) at any iteration by using the following proposition. 

\begin{proposition} \label{p:ValidInequalities} 
	Define $ \bar \varOmega^i= \left\{ S \subseteq I \setminus \{0\}:\ i \in S,\ |S| \leq \bar c \right\} $ and $ \underline \varOmega^i= \left\{ S \subseteq I \setminus \{0\}:\ i \notin S,\ |S| \leq \bar c \right\} $ as the set of feasible assortments that include and exclude product $ i \in I \setminus \{0\}$, respectively. Moreover, define $ \bar \pi^i=\max_{S \in \bar \varOmega^i} \pi(S) $ and $ \underline \pi^i=\max_{S \in \underline \varOmega^i} \pi(S) $. Then, inequalities $ \pi' \leq \bar \pi^i+1-x_i $ and $ \pi' \leq \underline \pi^i+x_i $ are valid to Problem $ (\ref{equ:equ32}) $ for all $ i \in I \setminus \{0\} $.
\end{proposition}

\proof{Proof.}
Suppose that assortment $ S $ does not include product $ i $, which means that $ x_i = 0 $. Then we have $ \pi (S) \leq \bar \pi^i + 1 $ and $ \pi (S) \leq \underline \pi^i $, where the first inequality holds because $ \pi(S) \leq 1 $ and $ \bar \pi^i \geq 0 $ and the second one holds because $ \pi(S) \leq \max_{S \in \underline \varOmega^i} \pi(S) $. If assortment $ S $ contains product $ i $, which means $ x_i = 1 $, then $ \pi (S) \leq \bar \pi^i $ and $ \pi (S) \leq \underline \pi^i + 1 $, where the first inequality holds because $ \pi(S) \leq \max_{S \in \bar \varOmega^i} \pi(S) $ and the second one holds because $ \pi(S) \leq 1 $ and $ \underline \pi^i \geq 0 $. $ \square $

In order to solve the bi-level problems required to calculate $ \bar \pi^i $ and $ \underline \pi^i $ for each $ i \in I \setminus \{0\} $, we use the dual Problem (\ref{equ:equ24}), which is a maximization problem so it can be used to transform $\max_{S \in \bar \varOmega^i} \pi(S) $ and $ \max_{S \in \underline \varOmega^i} \pi(S) $ into single-level \MIP{}s.

If $ \bar u = 1 $ in Problem (\ref{equ:equ32}), we can further strengthen the formulation by removing some of the products in $ I \setminus \{0\} $ that are not part of any optimal assortment. By doing so, we can significantly reduce the solution space for the next iterations. In this strategy, we use Algorithm \ref{alg:1} to find the set $\Lambda$ in polynomial time. Following the results from Section \ref{subsec:UnconSolution}, we know that an optimal assortment for the unconstrained version of the \APP{} with $ \bar u = 1 $ is in a set $ S'(\lambda)$, for some $\lambda \in \Lambda$. The following proposition describes the required conditions to permanently remove candidate products at intermediate iterations of Algorithm \ref{alg:2} when $ \bar u = 1 $. We define $ \underline z $ as a lower bound on the optimal objective value of Problem (\ref{equ:equ4}), which we initialize with $ \underline z = 0 $. We update the value of $ \underline z $ as we find assortments improving the revenue along iterations of Algorithm \ref{alg:2}. That is, if $ \underline z < \pi(S^\ast_T) \sum_{i \in S^\ast_T} \rho_i(S^\ast_T)r_i $ at any iteration $ T $, then we set $ \underline z = \pi(S^\ast_T) \sum_{i \in S^\ast_T} \rho_i(S^\ast_T)r_i $.

\begin{proposition} \label{p:SuperValid}
	Let $ \lambda_l $ be the largest value in $ \Lambda$ such that $ \lambda_l \leq \underline z $. Then, the optimal assortment to Problem $ (\ref{equ:equ4}) $ when $ \bar u = 1 $ excludes products in $ I \setminus (S'(\lambda_{l}) \cup \{0\})$ so they can be removed from the set of candidate products.
\end{proposition}

\proof{Proof.}
The result follows immediately by observing that if it is not optimal to include a product in the assortment in the unconstrained case, then it is also not optimal to do so in the cardinality-constrained case. Otherwise, the optimality of the unconstrained solution from Lemma \ref{l:lambda00} is contradicted. $\square$

\subsubsection{Greedy approach}$  $\\
\label{subsubsec:Greedy} 
The proposed solution approach for the constrained \APP{} relies on an \MIP{} formulation, thus it may not be possible to solve large-scale instances to optimality in reasonable time. For this reason, this section describes a greedy algorithm for the constrained version of \APP{} for a general value of $ \bar u $. Our computational results indicate that our greedy algorithm reaches a solution with a small gap relative to the optimal solution in reasonable time. Algorithm \ref{alg:greedy} describes our greedy approach.

In Line 1, we initialize a lower bound $ \underline z $ on the optimal assortment's profit. The loop in Lines $2$--$13$ constructs an assortment denoted by $ S_i $, for each $ i \in I \setminus \{0\} $. Lines $ 3$--$5 $ update the lower-bound profit $ \underline z $ by finding $ \pi(S_i) $ and calculating the profit associated with assortment $ S_i $. The while-loop in Lines $ 6$--$13$ adds the most profitable product until the assortment reaches its maximum size, $ \bar c $, or when there is no profitable product to add. Line $ 7 $ resets $ \hat k $, while Lines $ 8$-$12 $ evaluate the profit of adding each product not currently in the assortment one at a time. If adding a product improves the lower bound $ \underline z  $, Line 12 updates the value of $ \underline z $ and stores the product index in $ \hat k $. After reviewing all the products not currently in the assortment, Line $ 13 $ adds product $ \hat k $, if any, to the current assortment. If the algorithm visits Line $ 13 $ with $ \hat k = \texttt{NIL} $ (nonexistent), then no product is added and assortment $ S_i $ is finalized. Line $ 14 $ finds the assortment with the largest profit among $ S_1,...,S_{|I|-1} $. Line $ 15 $ returns the assortment with the largest profit. Algorithm \ref{alg:greedy} finishes in finite time and its worst-case complexity is $ O(|I|^3 \sigma^\pi(I)) $, where $ \sigma^\pi (S) $ is the complexity of finding $ \pi(S) $ using linear program (\ref{equ:equ24}) in Lines 4 and 9.   
 
\begin{algorithm} 
	\small{
		\caption{: Greedy algorithm for the constrained \APP{} \label{alg:greedy}}
		\begin{algorithmic}[1]
			\State Initialize $\underline z = 0$
			\ForAll{$i \in I \setminus \{0\}$}
			\State Initialize $S_i=\{i\}$ and $ \hat k = i $
			\State Use linear program $ (\ref{equ:equ24}) $ to find $ \theta(S_i) $ and set $ \pi(S_i) = e^ {\theta(S_i)} $        
			\State Set $ \underline z= \pi(S_i) \sum_{j \in S_i} \rho_j(S_i)r_j $
			\While {$ |S_i| \leq \bar c $ and $ \hat k \neq \texttt{NIL} $}
			\State Set $ \hat k = \texttt{NIL} $
			\ForAll{$k \in I \setminus (S_i \cup \{0\})$}
			\State Use linear program $ (\ref{equ:equ24}) $ to find $ \theta(S_i \cup \{k\}) $ and set $ \pi(S_i \cup \{k\}) = e^ {\theta(S_i \cup \{k\})} $
			\State Compute $ z = \pi(S_i \cup \{k\}) \sum_{j \in S_i \cup \{k\}} \rho_j(S_i \cup \{k\})r_j $
			\If {$ z>\underline z $}
			\State Set $ \underline z=z $ and $ \hat k = k $
			\EndIf
			\EndFor
			\State Set $ S_i = S_i \cup \{\hat k\} $
			\EndWhile
			\EndFor
			\State Set $ i^\ast = \argmax_{i \in I \setminus \{0\}} \pi(S_i) \sum_{j \in S_i} \rho_j(S_i)r_j $\\
			\Return $ S_{i^\ast} $
	\end{algorithmic}}
\end{algorithm}
\normalsize  

\subsubsection{Solving extensions of \APP{}} \label{subsubsec:MultiCustomerType} $ $

Although our exact \PLUB{} and greedy solution methods can be extended to solve many other variants of the \APP{}, in this section we focus on two of them that increase the realism in the assortment planning process. The first variant considers space limitations (e.g., shelf space). In this case, each product occupies a given volume or area and the retail store has finite space for the products \citep{gallego2014constrained}. To accommodate this constraint, we define $ \mu $ as the total available space and $ \gamma_i < \mu$ as the amount of space required by product $ i $. Using these elements, we modify the feasible region $X$ in Problem (\ref{equ:equ32}) by replacing the cardinality constraint with the constraint $\sum_{i=1}^n \gamma_i x_i \leq \mu $ and update sets $ \bar \varOmega^i= \left\{ S \subseteq I \setminus \{0\}:\ i \in S,\ \sum_{j \in S} \gamma_j \leq \mu \right\} $ and $ \underline \varOmega^i= \left\{ S \subseteq I \setminus \{0\}:\ i \notin S,\ \sum_{j \in S} \gamma_j \leq \mu \right\} $. Algorithm \ref{alg:2} remains unchanged. To modify Algorithm \ref{alg:greedy}, we remove the condition $|S_i| \leq \bar{c}$ from Line 6 and update the condition in Line $11$ to ``\texttt{\textbf{if} $z > \underline{z}$ and $\sum_{j \in S_i \cup \{k\}} \gamma_j \leq \mu$}''.

Our assortment problem can also be extended to include different customer categories (i.e., different market segments) with different choice behaviors for the same set of products. To this end, we define $ C $ as the set of customer categories and $ \bar w_c $ as the proportion of customers in category $ c \in C $.  We also define $ \bar u_c $, $ p_{ic}^{0} $, $ \eta_{ikc} $, $ \nu_{ic} $, $ \rho_{ic}(S) $,  and $ y_{ikc} $ as the $ \bar u$-, $ p^{0} $-, $ \eta $-, $ \nu $-, $ \rho(S) $-parameters and $ y$-decision variables for customer category $ c \in C $, respectively. Moreover, we define $Y_c(S)$ as the feasible set of worst-case preference lists for customer $c \in C$ given assortment $S$. Using this notation, Problem $ (\ref{equ:equ4}) $ can be reformulated as

\begin{equation}\label{equ:equ50}
\begin{aligned}
z^\ast= & \max _{S \in \varOmega} \left\{\sum_{c \in C} \bar w_c \left( \min_{\mathbf{y_c} \in Y_c(S)} \prod_{k=1}^{{\bar u_c}} \left( 1-\sum_{i \in \bar S} {\eta_{ikc} p_{ic}^{0} y_{ikc}} \right) \right) \sum_{i \in S} \rho_{ic}(S)r_i\right\}.
\end{aligned}
\end{equation}

Problem (\ref{equ:equ50}) allows the retailer to include category-specific parameters to model the customer tolerance when facing unavailable products, resulting in a worst-case preference order for each category. Although this problem is more realistic, it is computationally harder compared to that with only one customer type (see \citet{bront2009column}, \citet{feldman2015bounding}, and \citet{bergman2017discrete}). To solve Problem $ (\ref{equ:equ50}) $, we extend our \PLUB{} and greedy algorithms by defining $ \pi_c(S) $ as the value of $ \pi(S) $ for customer category $ c \in C$. For each customer category $ c \in C $, we also define the dual feasible region $A_c(\mathbf{x})$ and the piecewise linear upper bound $ \pi^\prime_{c} $ for $ \pi_c(S) $ when using the feasible set $ Q_T^c $ at iteration $ T $ of the \PLUB{} approach. As a result, the upper bound problem solved at iteration $ T $ of Algorithm \ref{alg:2} is defined as
%\begin{equation}\label{equ:equ51}
\begin{multline} \label{equ:equ51}
\bar z_T = \max \left\{ \sum_{c \in C} \bar w_c \delta_c: \sum_{j =1}^n \left({\nu_{jc} x_j \left(r_j \pi'_{c}-\delta_c \right)} \right)\geq \delta_c,\ \mathbf{x} \in X,\ \mathbf{\alpha_c}, \mathbf{\beta_c} \in A_c(\mathbf{x}),\ \delta_c \in \mathbb{R}, \right.\\ \biggl. (\mathbf{q}_c, \mathbf{w}_c, \pi'_{c}) \in Q^c_{T},\ \forall c \in C \biggr\}.
\end{multline}
%\end{equation}

Using the same linearization techniques as in Problem (\ref{equ:equ32}), we can transform Problem (\ref{equ:equ51}) into an \MIP{}. Moreover, the strengthening procedure introduced in Proposition \ref{p:ValidInequalities} can also be used in Problem $ (\ref{equ:equ51}) $ by replacing $ \bar \pi^i $ and $ \underline \pi^i $ for each customer category $ c $ with $ \bar \pi^i_c = \max_{S \in \bar \varOmega^i} \pi_c(S) $ and $ \underline \pi^i_c = \max_{S \in \underline \varOmega^i} \pi_c(S) $, respectively. In order to use our greedy approach in Algorithm \ref{alg:greedy} we calculate $ \underline z $ and $ z $ in Lines 5 and 10 using the objective function of Problem (\ref{equ:equ51}). To do so, we compute $ \pi(\cdot) $ and $ \theta(\cdot) $ for every customer category $ c \in C $, which we denote by $ \pi_c(\cdot) $ and $ \theta_c(\cdot) $, respectively. We provide an updated version of Algorithm \ref{alg:2} and \ref{alg:greedy} for the multi-category problem in the supplementary material.

\section{Computational experiments} \label{sec:Results}
In this section, we describe the procedure used to generate random instances and show the performance of our solution approaches when solving the cardinality-constrained \APP{}. Next, we analyze the impact that the customers' sensitivity to product unavailability has on the optimal assortment.    

\subsection{Random instance generation} \label{InstanceGeneration}

We generate problem instances of different sizes, combining values $ n \in \{20,40,60,80,100\} $, $ \bar c \in \{0.1n,0.2n,0.3n\} $, and $ \bar u \in \{1,2,3,4,5\} $. We generate 5 instances for every combination of $ n $, $ \bar c $, and $ \bar u $, resulting in 375 instances in total. To generate the remaining parameters, we first generate auxiliary parameters $ a_i $, $ b_i $, and $ d_i $ from a uniform distribution in $ [0.75,1.25] $ and parameter $ o_i $ from a uniform distribution in $ [0,1] $, for each product $ i \in I \setminus \{0\} $. Following the approach in \citet{gallego2014constrained}, we calculate the revenue and utility of product $ i $ as $ r_i=10 \times o_i^2 \times a_i $ and $ \nu_i=10 \times (1-o_i) \times b_i $, respectively. By using parameter $ o_i $ in both cases, we incorporate the documented negative correlation between revenue and utility of each product and also introduce some noise due to parameters $ a_i $ and $ b_i $ since such negative correlation does not always exist \citep{gallego2014constrained}. Because customers may immediately leave the store when a highly preferred product is unavailable, we generate leaving probabilities as $ p_i^0=0.4 \times (1-o_i) \times d_i $ for each $ i \in I \setminus \{0\} $. In order to generate $ \eta_{ik} $, we use the function
\begin{align*}
\eta_{ik}=\frac{2}{1+e^{-(k-1)(1-o_i)}},
\end{align*}
which is increasing in $k$, the number of visited unavailable products, and has a positive correlation with the leaving probability $ p_i^0 $ and utility $\nu_i$ due to the term $ (1-o_i)$. The mechanism to generate $p_i^0$ and $\eta_{ik}$ also guarantee that the leaving probability $p_i^0 \eta_{ik}$ is always in the interval $[0,1]$ for every product $i$ and position $k$ in the preference list.

Our algorithms are implemented in C++, using ILOG CPLEX Optimization Studio 12.7.1 as the solver on an Intel CPU Core i7-6700 with 3.4 GHz, 16 GB of RAM with one thread. We set a time limit of 3600 seconds for solving each instance, which includes all the steps required by our algorithms.

\subsection{Performance of the optimization algorithms} \label{AlgorithmPerformance}
Table \ref{time table} shows the average run time of our iterative \PLUB{} approach from Section $ \ref{subsec:ConSolution} $ for instances solved to optimality within the time limit. Table \ref{time table} also shows the algorithm performance when using the valid inequalities from Proposition \ref{p:ValidInequalities}, which we denote as \PLUB{}$+ \overline {\underline {\pi}} $. Using the \PLUB{} approach, we can solve 317 instances out of 375 (i.e., $ 84.5\% $) within the time limit. Using both the \PLUB{} and $ \overline {\underline \pi} $ bounds, this number increases to 330 (i.e., $ 88\% $). However, for small instances, such bounds seem to be less effective, and even slow down the solution time compared to \PLUB{}. The reason is the trade-off between strengthening the formulation with the $ \overline {\underline \pi} $-bounds and the time required to calculate them. Table \ref{time table} also shows that as the values of $ \bar u $ and $\bar{c}$ increase, less instances can be solved within the time limit, as expected. For smaller values of $ \bar u $, which is common for many types of products, we can solve all instances in reasonable time even for 100 alternatives and assortment cardinality limit of $ 30 $ products. 

\begin{table}[]
\centering
\small
\caption{Average run time (in seconds) using \PLUB{} from Algorithm \ref{alg:2} and \PLUB{}+$ \overline {\underline \pi} $ from Proposition \ref{p:ValidInequalities}}
	\label{time table}
\begin{tabular}{ccrlrllrlrllrlrl}
\hline
		\multicolumn{1}{l}{}  & \multicolumn{1}{l}{}       & \multicolumn{4}{c}{$ \bar c =0.1 n$}                                               &  & \multicolumn{4}{c}{$ \bar c =0.2 n$}                           &  & \multicolumn{4}{c}{$ \bar c =0.3 n$}                            \\ \cline{3-6} \cline{8-11} \cline{13-16} 
		\multicolumn{1}{l}{$ n $} & \multicolumn{1}{l}{$ \bar u $} & \multicolumn{2}{c}{\PLUB{}}                       & \multicolumn{2}{c}{\PLUB{}+$ \overline {\underline \pi} $} &  & \multicolumn{2}{c}{\PLUB{}}   & \multicolumn{2}{c}{\PLUB{}+$ \overline {\underline \pi} $} &  & \multicolumn{2}{c}{\PLUB{}}   & \multicolumn{2}{c}{\PLUB{}+$ \overline {\underline \pi} $} \\ \hline
20  & 1                     & 0.05            & (5)         & 0.21                        & (5)  &  & 0.07            & (5)         & 0.25                        & (5)  &  & 0.12            & (5)         & 0.31                      & (5)    \\
    & 2                     & 0.05            & (5)         & 0.31                        & (5)  &  & 0.07            & (5)         & 0.37                        & (5)  &  & 0.10            & (5)         & 0.50                      & (5)    \\
    & 3                     & 0.07            & (5)         & 0.40                        & (5)  &  & 0.11            & (5)         & 0.59                        & (5)  &  & 0.28            & (5)         & 0.78                      & (5)    \\
    & 4                     & 0.07            & (5)         & 0.53                        & (5)  &  & 0.26            & (5)         & 0.77                        & (5)  &  & 0.97            & (5)         & 1.33                      & (5)    \\
    & 5                     & 0.13            & (5)         & 0.60                        & (5)  &  & 0.42            & (5)         & 1.09                        & (5)  &  & 1.89            & (5)         & 2.18                      & (5)    \\
    \hline
40  & 1                     & 0.07            & (5)         & 0.70                        & (5)  &  & 0.23            & (5)         & 1.05                        & (5)  &  & 0.36            & (5)         & 1.27                      & (5)    \\
    & 2                     & 0.10            & (5)         & 1.26                        & (5)  &  & 0.46            & (5)         & 1.72                        & (5)  &  & 7.72            & (5)         & 5.87                      & (5)    \\
    & 3                     & 0.42            & (5)         & 1.92                        & (5)  &  & 5.43            & (5)         & 10.93                       & (5)  &  & 133.28          & (5)         & 71.80                     & (5)    \\
    & 4                     & 2.22            & (5)         & 3.38                        & (5)  &  & 137.62          & (5)         & 97.99                       & (5)  &  & 231.99          & (5)         & 44.72                     & (5)    \\
    & 5                     & 3.27            & (5)         & 4.84                        & (5)  &  & 67.73           & (5)         & 86.26                       & (5)  &  & 613.08          & (5)         & 404.91                    & (5)    \\
    \hline
60  & 1                     & 0.12            & (5)         & 1.50                        & (5)  &  & 0.43            & (5)         & 2.12                        & (5)  &  & 1.50            & (5)         & 2.83                      & (5)    \\
    & 2                     & 0.24            & (5)         & 3.41                        & (5)  &  & 1.38            & (5)         & 5.12                        & (5)  &  & 9.49            & (5)         & 8.93                      & (5)    \\
    & 3                     & 6.01            & (5)         & 7.03                        & (5)  &  & 422.14          & (5)         & 232.99                      & (5)  &  & 914.40          & (5)         & 745.92                    & (5)    \\
    & 4                     & 32.66           & (5)         & 43.22                       & (5)  &  & 1126.14         & (5)         & 478.24                      & (5)  &  & 1829.12         & (2)         & 1106.90                   & (3)    \\
    & 5                     & 115.11          & (5)         & 50.03                       & (5)  &  & 145.03          & (1)         & 778.10                      & (4)  &  & \_              & (0)         & 797.72                    & (1)    \\
    \hline
80  & 1                     & 0.49            & (5)         & 3.55                        & (5)  &  & 0.97            & (4)         & 4.69                        & (5)  &  & 17.14           & (5)         & 7.89                      & (5)    \\
    & 2                     & 1.57            & (5)         & 6.83                        & (5)  &  & 2.75            & (5)         & 12.81                       & (5)  &  & 388.29          & (4)         & 34.62                     & (5)    \\
    & 3                     & 46.58           & (5)         & 15.59                       & (5)  &  & 5.98            & (4)         & 27.35                       & (4)  &  & 14.20           & (4)         & 699.17                    & (5)    \\
    & 4                     & 3.47            & (4)         & 16.37                       & (4)  &  & 7.55            & (3)         & 65.11                       & (4)  &  & 60.35           & (4)         & 154.99                    & (4)    \\
    & 5                     & 10.00           & (4)         & 328.06                      & (4)  &  & 53.16           & (3)         & 86.08                       & (4)  &  & 16.30           & (4)         & 1435.31                   & (5)    \\
    \hline
100 & 1                     & 0.39            & (5)         & 5.07                        & (5)  &  & 1.05            & (5)         & 8.66                        & (5)  &  & 6.52            & (5)         & 15.34                     & (5)    \\
    & 2                     & 1.31            & (5)         & 13.63                       & (5)  &  & 7.22            & (5)         & 35.04                       & (5)  &  & 15.56           & (4)         & 91.29                     & (5)    \\
    & 3                     & 751.53          & (5)         & 86.76                       & (5)  &  & 166.81          & (2)         & \multicolumn{1}{l}{334.52}  & (2)  &  & 424.08          & (2)         & 995.72                    & (1)    \\
    & 4                     & 834.04          & (2)         & 650.51                      & (3)  &  & \_              & (0)         & \multicolumn{1}{l}{856.58}  & (1)  &  & \_              & (0)         & \_                        & (0)    \\
    & 5                     & 538.20          & (1)         & \multicolumn{1}{l}{278.05}  & (1)  &  & \_              & (0)         & \multicolumn{1}{l}{\_}      & (0)  &  & \_              & (0)         & \multicolumn{1}{r}{\_}    & (0)    \\ \hline
\end{tabular}
    
    \begin{tablenotes}
		\footnotesize
		\item \textit{Note:} Number of instances solved to optimality within 3600 seconds in parenthesis. Average run times calculated only for instances solved to optimality.
	\end{tablenotes}
\end{table}

Table $ \ref{greedy} $ shows the average run time of our greedy algorithm. In this table, we report the optimality gap only for those instances for which we obtain an optimal solution in Table \ref{time table} (either using the \PLUB{} or \PLUB{}+$ \overline {\underline \pi} $). The number of instances used to calculate the average gap is in parenthesis. On average, the greedy algorithm finds an optimal (or near optimal) solution with a gap that is no more than $ 4.68 \% $ within the time limit. The greedy algorithm reaches the optimal solution for 311 instances out of 332 instances for which we have an optimal solution.

\begin{table}[]
\centering
\footnotesize
\caption{Average run time (in seconds) and average optimality gap for Algorithm \ref{alg:greedy} (greedy)}
	\label{greedy}
\begin{tabular}{ccrlclrlclrlc}
		\hline
		\multicolumn{1}{l}{}  & \multicolumn{1}{l}{}       & \multicolumn{3}{c}{$ \bar c =0.1 n$}                                         &  & \multicolumn{3}{c}{$ \bar c =0.2 n$}                                         &  & \multicolumn{3}{c}{$ \bar c =0.3 n$}                                         \\ \cline{3-5} \cline{7-9} \cline{11-13} 
		\multicolumn{1}{l}{$ n $} & \multicolumn{1}{l}{$ \bar u $} & \multicolumn{2}{l}{Gap} & \multicolumn{1}{l}{Time (sec.)} &  & \multicolumn{2}{l}{Gap} & \multicolumn{1}{l}{Time (sec.)} &  & \multicolumn{2}{l}{Gap} & \multicolumn{1}{l}{Time (sec.)} \\ \hline
		20                    & 1                          & 0.00\%                   & (5)     & 0.13                           &  & 0.00\%                   & (5)     & 0.35                           &  & 0.00\%                   & (5)     & 0.54                           \\
		& 2                          & 0.00\%                   & (5)     & 0.16                           &  & 0.00\%                   & (5)     & 0.44                           &  & 0.00\%                   & (5)     & 0.69                           \\
		& 3                          & 0.00\%                   & (5)     & 0.18                           &  & 0.00\%                   & (5)     & 0.50                           &  & 0.00\%                   & (5)     & 0.79                           \\
		& 4                          & 0.00\%                   & (5)     & 0.21                           &  & 0.00\%                   & (5)     & 0.55                           &  & 0.00\%                   & (5)     & 0.87                           \\
		& 5                          & 0.00\%                   & (5)     & 0.30                           &  & 0.00\%                   & (5)     & 0.72                           &  & 0.89\%                   & (5)     & 1.07                           \\
		\hline
		40                    & 1                          & 0.00\%                   & (5)     & 1.96                           &  & 0.00\%                   & (5)     & 4.79                           &  & 0.00\%                   & (5)     & 7.20                           \\
		& 2                          & 0.00\%                   & (5)     & 2.62                           &  & 0.00\%                   & (5)     & 6.14                           &  & 0.00\%                   & (5)     & 9.41                           \\
		& 3                          & 0.00\%                   & (5)     & 3.12                           &  & 0.00\%                   & (5)     & 7.37                           &  & 0.29\%                   & (5)     & 11.51                          \\
		& 4                          & 0.00\%                   & (5)     & 3.68                           &  & 0.05\%                   & (5)     & 8.56                           &  & 4.68\%                   & (5)     & 13.14                          \\
		& 5                          & 0.00\%                   & (5)     & 4.28                           &  & 2.83\%                   & (5)     & 10.12                          &  & 3.04\%                   & (5)     & 14.89                          \\
		\hline
		60                    & 1                          & 0.00\%                   & (5)     & 11.76                          &  & 0.00\%                   & (5)     & 25.04                          &  & 0.00\%                   & (5)     & 37.99                          \\
		& 2                          & 0.00\%                   & (5)     & 16.00                          &  & 0.00\%                   & (5)     & 35.21                          &  & 0.00\%                   & (5)     & 54.48                          \\
		& 3                          & 0.00\%                   & (5)     & 18.94                          &  & 0.00\%                   & (5)     & 40.93                          &  & 1.50\%                   & (5)     & 62.55                          \\
		& 4                          & 0.00\%                   & (5)     & 22.72                          &  & 1.24\%                   & (5)     & 49.20                          &  & 0.00\%                   & (3)     & 74.68                          \\
		& 5                          & 0.28\%                   & (5)     & 26.88                          &  & 0.08\%                   & (4)     & 58.16                          &  & 0.00\%                   & (1)     & 86.84                          \\
		\hline
		80                    & 1                          & 0.14\%                   & (5)     & 37.17                          &  & 0.31\%                   & (5)     & 72.74                          &  & 0.54\%                   & (5)     & 112.85                         \\
		& 2                          & 0.80\%                   & (5)     & 49.20                          &  & 1.50\%                   & (5)     & 99.29                          &  & 0.00\%                   & (5)     & 150.33                         \\
		& 3                          & 0.00\%                   & (5)     & 63.26                          &  & 0.00\%                   & (4)     & 130.13                         &  & 0.20\%                   & (5)     & 195.09                         \\
		& 4                          & 0.00\%                   & (4)     & 79.58                          &  & 0.00\%                   & (4)     & 163.00                         &  & 0.00\%                   & (4)     & 239.73                         \\
		& 5                          & 0.00\%                   & (5)     & 98.40                          &  & 0.00\%                   & (4)     & 201.78                         &  & 0.00\%                   & (5)     & 296.32                         \\
		\hline
		100                   & 1                          & 0.00\%                   & (5)     & 94.92                          &  & 0.00\%                   & (5)     & 163.98                         &  & 0.00\%                   & (5)     & 238.27                         \\
		& 2                          & 0.00\%                   & (5)     & 123.51                         &  & 0.00\%                   & (5)     & 260.65                         &  & 0.00\%                   & (5)     & 344.35                         \\
		& 3                          & 0.00\%                   & (5)     & 161.61                         &  & 0.00\%                   & (2)     & 333.66                         &  & 0.00\%                   & (2)     & 447.21                         \\
		& 4                          & 0.00\%                   & (3)     & 211.91                         &  & 0.00\%                   & (1)     & 411.95                         &  & \_                       & (0)     & 570.69                         \\
		& 5                          & 0.00\%                   & (1)     & 252.34                         &  & \_                       & (0)     & 494.36                         &  & \_                       & (0)     & 720.51                         \\ \hline
	\end{tabular}
    
    \begin{tablenotes}
		\footnotesize
		\item \textit{Note:} Average run times calculated over 5 instances. Number of instances solved to optimality by \PLUB{} or \PLUB{}+$ \overline {\underline \pi} $ in parenthesis. Average gap calculated over instances with known optimal solution. The gap of each instance calculated as $100\times \frac{z^\ast-z^G}{z^\ast} $, where $z^G$ is the greedy objective value from Algorithm \ref{alg:greedy}.
	\end{tablenotes}
\end{table}

Considering the run time and solution quality of our greedy algorithm, we test whether using its best found assortment improves the performance of our exact \PLUB{} approach when solving hard instances (i.e., those with $n \in \{60,80,100\}$ and $\bar u \in \{3,4,5\}$). To achieve this, we use the best assortment found by Algorithm \ref{alg:greedy} (i.e., $ S_{i^\ast} $) as a starting assortment to tighten the upper bound Problem (\ref{equ:equ32}) in Algorithm \ref{alg:2}, where assortment $S_{i^\ast}$ is included in the construction of $Q_T$ for every iteration $T \geq 1$. We call this strategy G-\PLUB{}+$ \overline {\underline \pi} $, as we also use the valid inequalities from Proposition \ref{p:ValidInequalities}. Table \ref{InitializedWithGreedy} compares the performance of  G-\PLUB{}+$ \overline {\underline \pi} $ and \PLUB{}+$ \overline {\underline \pi} $ when attempting to solve the same instances in $3600$ seconds. In this case, G-\PLUB{}+$ \overline {\underline \pi} $ solves 3 instances that are not solved by \PLUB{}+$ \overline {\underline \pi} $ within the time limit. However, initializing Algorithm \ref{alg:2} with the greedy assortment deteriorates the average run time for some instances.

\begin{figure}[ht!]
\centering
\tiny
    \begin{subfigure}{0.49\textwidth}
    \centering
      \includegraphics[scale=0.45]{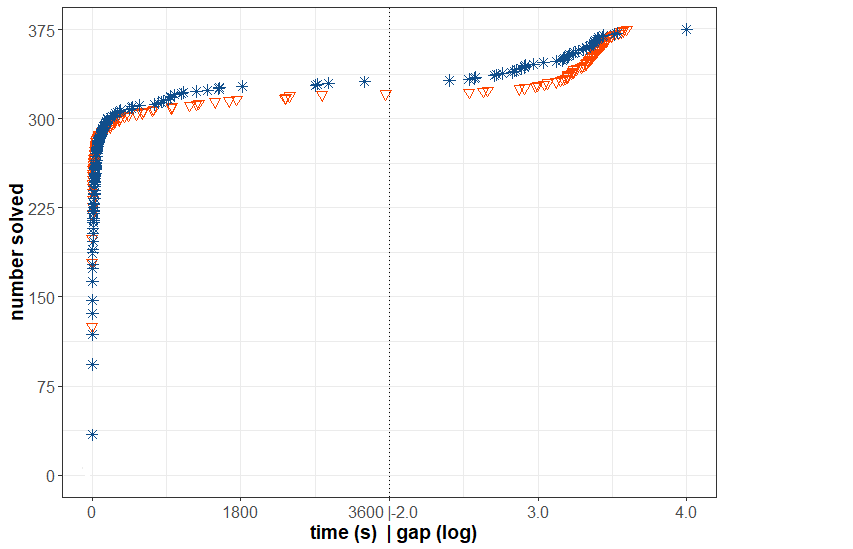}
      \caption{}
      \label{f:Plot-AllInstance}
    \end{subfigure}
    \hspace*{\fill}
    \begin{subfigure}{0.49\textwidth}
      \centering    
      \includegraphics[scale=0.38]{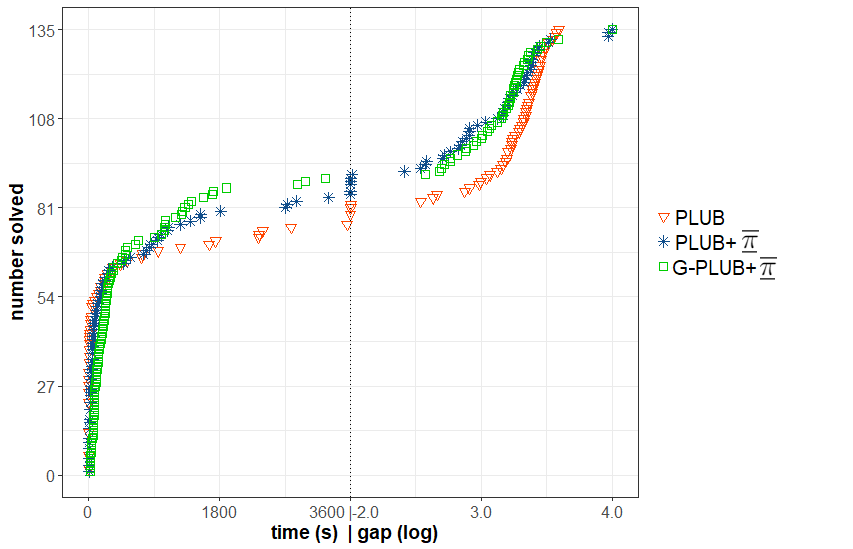}
      \caption{}
      \label{f:Plot-HardInstance}
    \end{subfigure}
    \caption{Cumulative performance plots to compare a) \PLUB{} and \PLUB{}+$ \overline{\underline \pi} $  over all instances and b) \PLUB{}, \PLUB{}+$ \overline{\underline \pi} $, and G-\PLUB{}+$ \overline{\underline \pi} $ over instances with $n \in \{60,80,100\}$ and $\bar u \in \{3,4,5\}$}
\end{figure}

We further explore the performance of \PLUB{}, \PLUB{}+$ \overline {\underline \pi} $, and G-\PLUB{}+$ \overline {\underline \pi} $ by analyzing the number of instances solved during the 3600 seconds limit. Figures \ref{f:Plot-AllInstance} and \ref{f:Plot-HardInstance} show cumulative performance plots for all 375 generated instances and for the subset of large instances with $n \in \{60,80,100\}$ and $\bar u \in \{3,4,5\}$, respectively. For each solution approach, the left half of each plot represents the number of instances solved within the time shown in the horizontal axis. The right half of each plot shows the number of instances that cannot be solved within 3600 seconds and that have $\log$ of the absolute gap (i.e., $\log_{10}(UB-LB)$) of at most the number shown on the horizontal axis. These figures show that \PLUB{}+$ \overline {\underline \pi} $ solves more instances to optimality after 900 seconds and also results in smaller optimality gaps than \PLUB{} for those instances not solved within the time limit. Figure (\ref{f:Plot-HardInstance}) also provides evidence of the superiority of \PLUB{}+$ \overline {\underline \pi} $ over \PLUB{} for large instances. Figure \ref{f:Plot-HardInstance} shows that G-\PLUB{}+$ \overline {\underline \pi} $ solves more instances to optimality after $900$ seconds compared to \PLUB{} and \PLUB{}+$ \overline {\underline \pi} $. Regarding the optimality gap after timing out, G-\PLUB{}+$ \overline {\underline \pi} $ and \PLUB{}+$ \overline {\underline \pi} $ have very similar performance. 

\begin{landscape}
\begin{table}[]
\centering
\vspace{50mm}
\small
\caption{Performance of \PLUB{}+$ \overline {\underline \pi} $ approach initialized with the greedy assortment from Algorithm \ref{alg:greedy}}
	\label{InitializedWithGreedy}
\begin{tabular}{ccrlcrlrccrlrlcrl}
\hline
    &                        & \multicolumn{5}{c}{$\bar c = 0.1n$}                                             & \multicolumn{5}{c}{$\bar c = 0.2n$}                                    & \multicolumn{5}{c}{$\bar c = 0.3n$}                                              \\ \cline{3-17} 
$n$   & $\bar u$  & \multicolumn{2}{c}{\PLUB{}+$ \overline {\underline \pi} $} & Greedy & \multicolumn{2}{c}{G-\PLUB{}+$ \overline {\underline \pi} $}                     & \multicolumn{2}{c}{\PLUB{}+$ \overline {\underline \pi} $} & Greedy & \multicolumn{2}{c}{G-\PLUB{}+$ \overline {\underline \pi} $}            & \multicolumn{2}{c}{\PLUB{}+$ \overline {\underline \pi} $}  & Greedy & \multicolumn{2}{c}{G-\PLUB{}+$ \overline {\underline \pi} $}                     \\ \hline
60  & \multicolumn{1}{c|}{3} & 7.03           & (5)          & 18.94  & 3.66   & \multicolumn{1}{l|}{(5)}          & 232.99     & (5)     & 40.93  & 236.88          & \multicolumn{1}{l|}{(5)} & 745.92          & (5)          & 62.55  & 542.27 & \multicolumn{1}{l|}{(5)}          \\
    & \multicolumn{1}{c|}{4} & 43.22          & (5)          & 22.72  & 26.05 & \multicolumn{1}{l|}{(5)}          & 478.24     & (5)     & 49.20  & 629.25          & \multicolumn{1}{l|}{(5)} & 1106.90         & (3)          & 74.56  & 1513.96         & \multicolumn{1}{l|}{(4)} \\
    & \multicolumn{1}{c|}{5} & 50.03 & (5)          & 26.88  & 91.11           & \multicolumn{1}{l|}{(5)}          & 778.10     & (4)     & 58.11  & 1312.46         & \multicolumn{1}{l|}{(4)} & 797.72 & (1)          & 85.88  & 1238.59         & \multicolumn{1}{l|}{(1)}          \\ \hline
80  & \multicolumn{1}{c|}{3} & 15.59          & (5)          & 63.26  & 2.47   & \multicolumn{1}{l|}{(5)}          & 27.35               & (4)     & 130.13 & 0.43   & \multicolumn{1}{l|}{(4)} & 699.17          & (5) & 194.29 & 11.23           & \multicolumn{1}{l|}{(4)}          \\
    & \multicolumn{1}{c|}{4} & 16.37          & (4)          & 79.83  & 0.75   & \multicolumn{1}{l|}{(4)}          & 65.11               & (4)     & 163.34 & 38.14  & \multicolumn{1}{l|}{(4)} & 154.99          & (4)          & 239.54 & 85.67  & \multicolumn{1}{l|}{(4)}          \\
    & \multicolumn{1}{c|}{5} & 328.06         & (4)          & 98.40  & 271.12 & \multicolumn{1}{l|}{(5)} & 86.08               & (4)     & 202.31 & 60.29  & \multicolumn{1}{l|}{(4)} & 1435.31         & (5) & 297.21 & 939.76          & \multicolumn{1}{l|}{(4)}          \\ \hline
100 & \multicolumn{1}{c|}{3} & 86.76          & (5)          & 161.61 & 50.15  & \multicolumn{1}{l|}{(5)}          & 334.52              & (2)     & 337.04 & 206.98 & \multicolumn{1}{l|}{(2)} & 995.72          & (1)          & 448.05 & 1813.95         & \multicolumn{1}{l|}{(2)} \\
    & \multicolumn{1}{c|}{4} & 650.51         & (3) & 213.92 & 221.47 & \multicolumn{1}{l|}{(2)}          & 856.58              & (1)     & 426.02 & 634.13 & \multicolumn{1}{l|}{(1)} & \_              & (0)          & \_     & \_              & \multicolumn{1}{l|}{(0)}          \\
    & \multicolumn{1}{c|}{5} & 278.05         & (1)          & 264.20 & 203.95 & \multicolumn{1}{l|}{(1)}          & \_                  & (0)     & \_     & \_              & \multicolumn{1}{l|}{(0)} & \_              & (0)          & \_     & \_              & \multicolumn{1}{l|}{(0)}          \\ \hline
\end{tabular}
    
    \begin{tablenotes}
		\footnotesize
		\item \textit{Note:} Number of instances solved to optimality within 3600 seconds by \PLUB{}+$ \overline {\underline \pi} $ and G-\PLUB{}+$ \overline {\underline \pi} $ in parenthesis. G-\PLUB{}+$ \overline {\underline \pi} $ column excludes the run time for Algorithm \ref{alg:greedy} (greedy). Average run times calculated only for instances solved to optimality.
	\end{tablenotes}
\end{table}
\end{landscape}

\subsection{Assortment analysis}

In this section, we analyze the variation in the optimal assortment due to changes in some parameters of the \APP{}. To measure such variation, we define $\mathcal{T}$ as a set of indices identifying a combination of parameters used. For every $t \in \mathcal{T}$, $S^\ast_t$ represents the corresponding optimal assortment when solved with parameters $t$. For two distinct indices $t,t^\prime \in \mathcal{T}$, $|S^\ast_t \setminus S^\ast_{t^\prime}| + |S^\ast_{t^\prime} \setminus S^\ast_t | $ calculates the number of products that are not in common in assortments $S^\ast_t$ and $S^\ast_{t^\prime}$. We  normalize $|S^\ast_t \setminus S^\ast_{t^\prime}| + |S^\ast_{t^\prime} \setminus S^\ast_t |$ by dividing it by $|S^\ast_t|+|S^\ast_{t^\prime}|$, its maximum value, and define an assortment variation index as
\begin{align*}
AV = \frac{1}{|\mathcal{T}|(|\mathcal{T}|-1)} \sum_{t,t^\prime \in \mathcal{T}} \frac{|S^\ast_t \setminus S^\ast_{t^\prime}| + |S^\ast_{t^\prime} \setminus S^\ast_t |}{|S^\ast_t|+|S^\ast_{t^\prime}|}.   
\end{align*}

We first assess how sensitive the optimal assortment is to the length of the top-priority preference list, $\bar u$, and the customers' level of sensitivity to product unavailability. To do so, we vary $\bar u$ from $ 1$ to $5$ for each instance, i.e., we construct five parameter combinations for the same instance, thus $|\mathcal{T}|=5$. We examine problem instances with $n=20$ and $\bar c \in \{2,6,10\}$ for $ 7 $ product categories (including $ 3 $ base and $ 4 $ mixed categories). Product Categories 1, 2, and 3 are base categories and consist of products for which customers are high-sensitive, medium-sensitive, and low-sensitive to product unavailability, respectively. Product categories 4--7 are mixed product categories, including products that significantly vary in price and/or quality, and therefore, customers have different sensitivity levels to their unavailability. We generate 5 instances for every combination of $n$, $\bar c$, and product category. We use the same instance generation process as in Section \ref{InstanceGeneration} with the exception of generating $ o_i $ for every product $i$ from a uniform distribution in $[0,0.2]$, $[0.4,0.6]$, and $[0.8,1]$, for base Categories 1, 2, and 3, respectively. Product Categories $4$, $5$, and $6$ consist of mixed Categories $\{1,2\}$, $\{1,3\}$, and $\{2,3\}$, respectively. For mixed Category $ \{1,2\} $, we generate $ o_i $ for every product $ i $ with probability $ 0.5 $ from a uniform distribution in $[0,0.2]$ and with probability $ 0.5 $ from a uniform distribution in $[0.4,0.6]$. Similarly, we generate $ o_i $ for products in Categories $ \{1,3\} $ and $ \{2,3\} $ from their corresponding ranges. The last category of products that we consider is a mixture of all base categories of products, for which we generate $ o_i $ for every product $ i $ with equal probability ($ 0.\bar3 $) using a uniform distribution in $ [0,0.2] $, $ [0.4,0.6] $, or $ [0.8,1] $.

Table \ref{AssortmentVariation} shows that $ AV $ increases with the number of products allowed in the assortment, meaning that an increase in $\bar{u}$ requires adjustments in the assortment to mitigate the impact that unavailable products have on the revenue, which can only be achieved when the cardinality constraint allows it. Moreover, as more products for which customers are less sensitive to their unavailability (e.g., Categories $3$ and $\{2,3\}$ compared to Categories $1$ and $\{1,2\}$) are included in the set of candidate products, the optimal assortment varies less with respect to $\bar u$, as expected. Table \ref{AssortmentVariation} also shows the average solution time for different product categories and $\bar c$. We see that, as expected, the run times increase as $\bar{c}$ increases.

\begin{table}[]
\centering
\small
\caption{Assortment variation ($AV$) with respect to $\bar u$ and run time (in seconds) }
	\label{AssortmentVariation}
\begin{tabular}{ccccccccc}
\hline
             & \multicolumn{2}{c}{$n=20$, $\bar c=2$} &  & \multicolumn{2}{c}{$n=20$, $\bar c=6$} &  & \multicolumn{2}{c}{$n=20$, $\bar c=10$} \\ \cline{2-3} \cline{5-6} \cline{8-9} 
Product Category & $AV$             & Time (sec.) &  & $AV$ & Time (sec.) &  & $AV$             & Time (sec.)  \\ \hline
1            & 0.24  & 0.11        &  & 0.42  & 0.50        &  & 0.51  & 2.54         \\
2            & 0     & 0.07        &  & 0.24  & 0.53        &  & 0.48  & 2.35         \\
3            & 0     & 0.07        &  & 0.16  & 0.28        &  & 0.30  & 0.79         \\ \hline
\{1,2\}          & 0.12  & 0.10        &  & 0.28  & 1.53        &  & 0.46  & 12.09        \\
\{1,3\}          & 0    & 0.05        &  & 0.16  & 0.58        &  & 0.44  & 1.26         \\
\{2,3\}          & 0     & 0.06        &  & 0.12  & 0.15        &  & 0.36  & 0.49         \\ \hline
\{1,2,3\}        & 0     & 0.69        &  & 0.22  & 1.10        &  & 0.60  & 1.27         \\ \hline
\end{tabular}
    
    \begin{tablenotes}
		\small
		\item \textit{Note:} Run times are averaged over 25 instances (5 instances per $\bar u \in \{1,...,5\}$). $AV$ values for each $\bar c$ and product category are averaged over 5 instances. Results are obtained by running \PLUB{}.
	\end{tablenotes}
\end{table}

Table \ref{AssortmentVariationMNL} compares the performance of \APP{} with the expected revenue maximization problem under the standard MNL model with cardinality constraint (AP-MNL), which is equivalent to our cardinality constrained \APP{} when $\bar u = 0$. We define $S_{MNL}$ as the optimal assortment to AP-MNL and $z_{MNL}$ as the corresponding objective value of $S_{MNL}$ in \APP{}. Table \ref{AssortmentVariationMNL} shows that for small values of $\bar c$ and for product categories including low-sensitive customers, AP-MNL is a good approximation to \APP{} given the small differences in the optimal assortments. However, for product categories with high-sensitive customers, using the $S_{MNL}$ can cause large revenue losses for the retailer as $\bar c$ increases.

\begin{table}[]
\centering
\small
\caption{Comparison between \APP{} and \textbf{AP-MNL}}
	\label{AssortmentVariationMNL}
\begin{tabular}{cclcclcclcc}
\hline
             &        &  & \multicolumn{2}{c}{$n=20,\bar c=2$} &  & \multicolumn{2}{c}{$n=20,\bar c=6$} &  & \multicolumn{2}{c}{$n=20,\bar c=10$} \\ \cline{4-5} \cline{7-8} \cline{10-11} 
Product Category & $\bar u$ &  & $AV$                        & Gap              &  & $AV$                       & Gap               &  & $AV$                        & Gap               \\ \hline
1            & 2      &  & 0.13             & 2.4\%            &  & 0.40            & 14.0\%            &  & 0.43             & 17.3\%            \\
             & 4      &  & 0.40             & 7.8\%            &  & 0.47            & 57.9\%            &  & 0.63             & 69.3\%            \\\hline
2            & 2      &  & 0                & 0\%              &  & 0.20            & 1.0\%             &  & 0.25             & 1.5\%             \\
             & 4      &  & 0                & 0\%              &  & 0.36            & 4.9\%             &  & 0.60             & 11.1\%            \\\hline 
3            & 2      &  & 0                & 0\%              &  & 0               & 0\%               &  & 0.02             & 0.2\%             \\
             & 4      &  & 0                & 0\%              &  & 0.20           & 1.0\%             &  & 0.39             & 6.1\%             \\\hline
\{1,2\}          & 2      &  & 0                & 0\%              &  & 0.34            & 5.5\%             &  & 0.40             & 10.7\%            \\
             & 4      &  & 0.20             & 9.8\%            &  & 0.35            & 47.0\%            &  &0.41             & 62.0\%            \\\hline
\{1,3\}          & 2      &  & 0                & 0\%              &  & 0.02            & 0.1\%             &  & 0.02             & 0.1\%             \\
             & 4      &  & 0                & 0\%              &  & 0.12            & 2.2\%             &  & 0.50             & 36.3\%            \\\hline
\{2,3\}          & 2      &  & 0                & 0\%              &  & 0               & 0\%               &  & 0                & 0\%               \\
             & 4      &  & 0                & 0\%              &  & 0.20            & 0.3\%             &  & 0.20             & 7.8\%             \\\hline
\{1,2,3\}        & 2      &  &0                & 0\%              &  & 0               & 0\%               &  & 0.02            & 0.02\%            \\
             & 4      &  & 0                & 0\%              &  & 0.33            & 14.7\%            &  & 1.00             & 59.9\%            \\ \hline
\end{tabular}
    
    \begin{tablenotes}
		\small
		\item \textit{Note:} $AV = \frac{|S^\ast \setminus S_{MNL}| + |S_{MNL} \setminus S^\ast |}{|S^\ast|+|S_{MNL}|}$ and $Gap = 100 \times \frac{z^\ast-z_{MNL}}{z^\ast}$, where $S^\ast$ is the optimal assortment for the constrained \APP{}. Shown $AV$ and $Gap$ are averaged over 5 instances.
	\end{tablenotes}
\end{table}

Table \ref{MutiCustomer} shows how the optimal assortment varies assuming a single customer category when in fact there are more than one (i.e., $|C| > 1$). Let $S^\ast$ be the optimal assortment to Problem (\ref{equ:equ50}) considering all categories of customers in $C$ with a given distribution $\bar w$ and define $S_{c}$ to denote the optimal assortment to Problem (\ref{equ:equ4}) assuming a single customer category. We define $AV = \frac{1}{|C|} \sum_{c \in C} \frac{|S^\ast \setminus S_{c}| + |S_{c} \setminus S^\ast |}{|S^\ast|+|S_{c}|}$ to show how $S_c$ varies from $S^\ast$ on average considering all $c \in C$. In Table \ref{MutiCustomer}, we use the previously defined high-, medium-, and low-sensitive customer categories (i.e., Categories 1, 2, and 3) and suppose that for every candidate product, customers are distributed into either two or all three of such categories. The optimal assortment $S^\ast$ changes as the customer distribution changes, sometimes very drastically, with respect to the case where customers are believed to belong to only one customer category. This is evidenced not only by the positive value of $AV$ but also by the positive (and sometimes high) value of the average gap, which indicates a revenue loss of more than 26\% in some cases when the existing multiple customer categories are ignored.

\begin{table}[]
\centering
\small
\caption{Assortment variation and run time (in seconds) when considering different customer categories}
	\label{MutiCustomer}
\begin{tabular}{cccccclccclccc}
\hline
\multirow{2}{*}{\begin{tabular}[c]{@{}c@{}}Customer\\ Categories ($C$)\end{tabular}} & \multirow{2}{*}{\begin{tabular}[c]{@{}c@{}}Distribution\\ ($\bar w$)\end{tabular}} & \multirow{2}{*}{$\bar u$} & \multicolumn{3}{c}{$n=20,\bar c = 2$} &  & \multicolumn{3}{c}{$n=20,\bar c = 6$} &  & \multicolumn{3}{c}{$n=20,\bar c = 10$} \\ \cline{4-6} \cline{8-10} \cline{12-14} 
                                                                             &                               &                                        & $AV$             & Gap               & Time          &  & $AV$             & Gap              & Time           &  & $AV$             & Gap              & Time           \\ \hline
1,2                                                                          & 0.2,0.8                       & 2                                      & 0.03           & 0.06\%            & 0.18          &  & 0.07           & 1.83\%           & 0.90           &  & 0.08           & 2.20\%           & 2.50           \\
                                                                             &                               & 4                                      & 0.18           & 16.99\%           & 0.24          &  & 0.25           & 16.66\%          & 2.40           &  & 0.31           & 22.79\%          & 6.97           \\ 
                                                                             & 0.5,0.5                       & 2                                      & 0.30           & 0\%               & 0.19          &  & 2.07           & 1.41\%           & 0.83           &  & 4.96           & 1.69\%           & 2.19           \\
                                                                             &                               & 4                                      & 0.29           & 16.58\%           & 0.22          &  & 3.88           & 15.24\%          & 2.57           &  & 10.98          & 20.69\%          & 8.24           \\
                                                                             & 0.8,0.2                       & 2                                      & 0.03           & 0.12\%            & 0.21          &  & 0.10           & 0.90\%           & 1.24           &  & 0.11           & 1.02\%           & 3.32           \\
                                                                             &                               & 4                                      & 0.18           & 15.07\%           & 0.25          &  & 0.27           & 10.81\%          & 3.89           &  & 0.35           & 14.41\%          & 18.70          \\ \hline
1,3                                                                          & 0.2,0.8                       & 2                                      & 0.22           & 6.50\%            & 0.17          &  & 0.18           & 6.16\%           & 0.45           &  & 0.19           & 6.71\%           & 0.65           \\
                                                                             &                               & 4                                      & 0.30           & 26.62\%           & 0.19          &  & 0.19           & 15.81\%          & 0.51           &  & 0.30           & 23.27\%          & 0.88           \\
                                                                             & 0.5,0.5                       & 2                                      & 0.22           & 4.87\%            & 0.16          &  & 0.18           & 4.76\%           & 0.59           &  & 0.20           & 5.21\%           & 1.18           \\
                                                                             &                               & 4                                      & 0.30           & 26.39\%           & 0.20          &  & 0.19           & 15.35\%          & 0.62           &  & 0.30           & 22.54\%          & 1.22           \\
                                                                             & 0.8,0.2                       & 2                                      & 0.22           & 2.91\%            & 0.20          &  & 0.18           & 2.41\%           & 1.02           &  & 0.20           & 2.65\%           & 2.99           \\
                                                                             &                               & 4                                      & 0.30           & 25.51\%           & 0.23          &  & 0.19           & 13.81\%          & 0.97           &  & 0.30           & 19.92\%          & 2.18           \\ \hline
2,3                                                                          & 0.2,0.8                       & 2                                      & 0.25           & 3.87\%            & 0.19          &  & 0.21           & 3.86\%           & 4.22           &  & 0.21           & 4.13\%           & 13.50          \\
                                                                             &                               & 4                                      & 0.25           & 4.68\%            & 0.31          &  & 0.23           & 5.52\%           & 5.67           &  & 0.24           & 8.63\%           & 14.34          \\
                                                                             & 0.5,0.5                       & 2                                      & 0.25           & 2.41\%            & 0.20          &  & 0.25           & 2.85\%           & 6.14           &  & 0.21           & 3.06\%           & 31.61          \\
                                                                             &                               & 4                                      & 0.25           & 3.59\%            & 0.31          &  & 0.23           & 4.09\%           & 6.78           &  & 0.25           & 6.46\%           & 29.11          \\
                                                                             & 0.8,0.2                       & 2                                      & 0.25           & 1.07\%            & 0.21          &  & 0.23           & 2.09\%           & 8.18           &  & 0.21           & 2.24\%           & 45.92          \\
                                                                             &                               & 4                                      & 0.25           & 1.90\%            & 0.39          &  & 0.24           & 2.23\%           & 14.68          &  & 0.27           & 4.17\%           & 81.33          \\ \hline
1,2,3                                                                        & 0.1,0.3,0.6                   & 2                                      & 0.02           & 0.05\%            & 0.44          &  & 0.04           & 1.29\%           & 1.52           &  & 0.05           & 1.55\%           & 3.53           \\
                                                                             &                               & 4                                      & 0.12           & 11.44\%           & 0.58          &  & 0.19           & 11.66\%          & 3.72           &  & 0.28           & 16.16\%          & 10.04          \\
                                                                             & $0.\bar3,0.\bar3,0.\bar3$                   & 2                                      & 0.02           & 0.03\%            & 0.49          &  & 0.04           & 1.12\%           & 1.76           &  & 0.05           & 1.35\%           & 3.63           \\
                                                                             &                               & 4                                      & 0.12           & 11.32\%           & 0.60          &  & 0.19           & 11.16\%          & 3.78           &  & 0.28           & 15.27\%          & 10.78          \\
                                                                             & 0.6,0.3,0.1                   & 2                                      & 0.04           & 0.04\%            & 0.41          &  & 0.07           & 0.88\%           & 2.42           &  & 0.07           & 1.04\%           & 4.98           \\
                                                                             &                               & 4                                      & 0.12           & 10.96\%           & 0.59          &  & 0.19           & 9.85\%           & 5.41           &  & 0.26           & 13.45\%          & 15.71          \\ \hline
\end{tabular}
    
    \begin{tablenotes}
		\footnotesize
		\item \textit{Note:} 
		$Gap = \frac{1}{|C|} \sum_{c \in C} \left(100 \times \frac{z^\ast-z_{c}}{z^\ast} \right)$, where $z_{c}$ is the objective value in Problem (\ref{equ:equ50}) corresponding to assortment $S_{c}$. $AV$ values, gaps, and run times are averaged over 5 instances.
	\end{tablenotes}
\end{table}

\section{Conclusions} \label{sec:Conclusion}
We study an assortment planning problem for substitutable products under a customer ranking-based choice model. This choice model incorporates two factors supported by empirical literature. The first is the effect that unavailable products have on the customer purchasing decision. In particular, we show that the sequence in which the customer attempts to buy unavailable top-priority products (i.e., the sequence in which such products are in the customer's top-priority preference list) has a big influence in the probability of leaving the store without purchase and must be considered in the assortment decisions. The second effect is that customers attempt to purchase the products in their preference list first rather than checking all the offered products in the store, which they do only when their most preferred products are unavailable and they decide to stay in the store. 

Because of the existing uncertainty in the estimation of the customers' preference lists, we focus on the worst-case preference list corresponding to each given assortment. To find the optimal assortment, we present a new mathematical programming formulation for our choice model and embed it into a bi-level optimization approach that maximizes retailer's expected revenue while considering the worst-case substitution behavior of the customers. Considering the complexity of the problem, we propose a set of exact and greedy solution methods. We consider both unconstrained and constrained versions of the problem, where the latter can be a cardinality constraint limiting the number of products in the assortment or a knapsack constraint reflecting space limitations. Our solution approach for the cardinality-constrained problem transforms the bi-level problem into a single-level upper bound problem, which is iteratively improved using cutting planes. We show that using this approach, we can solve instances of different sizes in reasonable times. We also provide a greedy algorithm that can solve large instances in reasonable time with small optimality gaps.

From a management perspective, our models, algorithms, and computational experience suggest that the effect of an unavailable product cannot be considered in isolation because the customer's behavior is also influenced by which and how many other products in the customer's preference list  are available. This means that the sequence of purchase attempts must be considered in the assortment planning process, otherwise the customer's leaving probability could be drastically  underestimated. Moreover, the cumulative unavailability effect cannot be ignored in the assortment decision, as it may result in a very different optimal assortment and a loss of revenue compared to the case when customers are assumed to have a sequence-independent purchasing behavior.  Additionally, incorporating customer categories with heterogeneous behaviors results in a more realistic analysis that balances the tradeoff between the importance of each category  (i.e., percentage of customers) and their sensitivity to unavailable products, ultimately producing a more realistic assortment. As expected in any assortment planning model, an increase in the assortment budget (e.g., cardinality or knapsack) allows the retailer to better accommodate the customer preferences and hedge the risk of having customers leaving the store due to product unavailability.  

%\section*{Acknowledgments}
%The authors thank Scott Webster, Professor of Supply Chain Management at W. P. Carey School of Business at Arizona State University, whose comments and support improved the paper in many ways.

\SingleSpacedXI
\footnotesize
\bibliographystyle{informs2014}
\bibliography{main}

\begin{thebibliography}{41}
\providecommand{\natexlab}[1]{#1}
\providecommand{\url}[1]{\texttt{#1}}
\providecommand{\urlprefix}{URL }

\bibitem[{Anderson et~al.(2006)Anderson, Fitzsimons, \protect\BIBand{}
  Simester}]{anderson2006measuring}
Anderson ET, Fitzsimons GJ, Simester D (2006) Measuring and mitigating the
  costs of stockouts. \emph{Management Science} 52(11):1751--1763.

\bibitem[{Aouad et~al.(2018)Aouad, Farias, Levi, \protect\BIBand{}
  Segev}]{aouad2018approximability}
Aouad A, Farias V, Levi R, Segev D (2018) The approximability of assortment
  optimization under ranking preferences. \emph{Operations Research}
  66(6):1661--1669.

\bibitem[{Berbeglia(2016)}]{berbeglia2016discrete}
Berbeglia G (2016) Discrete choice models based on random walks.
  \emph{Operations Research Letters} 44(2):234--237.

\bibitem[{Bergman \protect\BIBand{} Cire(2017)}]{bergman2017discrete}
Bergman D, Cire AA (2017) Discrete nonlinear optimization by state-space
  decompositions. \emph{Management Science} 64(10):4700--4720.

\bibitem[{Bernstein et~al.(2015)Bernstein, K{\"o}k, \protect\BIBand{}
  Xie}]{bernstein2015dynamic}
Bernstein F, K{\"o}k AG, Xie L (2015) Dynamic assortment customization with
  limited inventories. \emph{Manufacturing \& Service Operations Management}
  17(4):538--553.

\bibitem[{Besbes \protect\BIBand{} Saur{\'e}(2016)}]{besbes2016product}
Besbes O, Saur{\'e} D (2016) Product assortment and price competition under
  multinomial logit demand. \emph{Production and Operations Management}
  25(1):114--127.

\bibitem[{Birbil et~al.(2009)Birbil, Frenk, Gromicho, \protect\BIBand{}
  Zhang}]{birbil2009role}
Birbil {\c{S}}{\.I}, Frenk J, Gromicho JA, Zhang S (2009) The role of robust
  optimization in single-leg airline revenue management. \emph{Management
  Science} 55(1):148--163.

\bibitem[{Blanchet et~al.(2016)Blanchet, Gallego, \protect\BIBand{}
  Goyal}]{blanchet2016markov}
Blanchet J, Gallego G, Goyal V (2016) A markov chain approximation to choice
  modeling. \emph{Operations Research} 64(4):886--905.

\bibitem[{Bront et~al.(2009)Bront, M{\'e}ndez-D{\'\i}az, \protect\BIBand{}
  Vulcano}]{bront2009column}
Bront JJM, M{\'e}ndez-D{\'\i}az I, Vulcano G (2009) A column generation
  algorithm for choice-based network revenue management. \emph{Operations
  Research} 57(3):769--784.

\bibitem[{Chen \protect\BIBand{} Chen(2018)}]{chen2018robust}
Chen M, Chen ZL (2018) Robust dynamic pricing with two substitutable products.
  \emph{Manufacturing \& Service Operations Management} 20(2):249--268.

\bibitem[{Davis et~al.(2014)Davis, Gallego, \protect\BIBand{}
  Topaloglu}]{davis2014assortment}
Davis JM, Gallego G, Topaloglu H (2014) Assortment optimization under variants
  of the nested logit model. \emph{Operations Research} 62(2):250--273.

\bibitem[{D{\'e}sir et~al.(2019)D{\'e}sir, Goyal, Jiang, Xie, \protect\BIBand{}
  Zhang}]{desir2019nonconvex}
D{\'e}sir A, Goyal V, Jiang B, Xie T, Zhang J (2019) A nonconvex min-max
  framework and its applications to robust assortment optimization
  \urlprefix\url{https://faculty.insead.edu/antoine-desir/documents/Operations-Research-RMC_20181001.pdf}.

\bibitem[{Farias et~al.(2013)Farias, Jagabathula, \protect\BIBand{}
  Shah}]{farias2013nonparametric}
Farias VF, Jagabathula S, Shah D (2013) A nonparametric approach to modeling
  choice with limited data. \emph{Management Science} 59(2):305--322.

\bibitem[{Feldman \protect\BIBand{} Topaloglu(2015)}]{feldman2015bounding}
Feldman J, Topaloglu H (2015) Bounding optimal expected revenues for assortment
  optimization under mixtures of multinomial logits. \emph{Production and
  Operations Management} 24(10):1598--1620.

\bibitem[{Fitzsimons(2000)}]{fitzsimons2000consumer}
Fitzsimons GJ (2000) Consumer response to stockouts. \emph{Journal of consumer
  research} 27(2):249--266.

\bibitem[{Flores et~al.(2019)Flores, Berbeglia, \protect\BIBand{}
  Van~Hentenryck}]{flores2019assortment}
Flores A, Berbeglia G, Van~Hentenryck P (2019) Assortment optimization under
  the sequential multinomial logit model. \emph{European Journal of Operational
  Research} 273(3):1052--1064.

\bibitem[{Gallego \protect\BIBand{} Topaloglu(2014)}]{gallego2014constrained}
Gallego G, Topaloglu H (2014) Constrained assortment optimization for the
  nested logit model. \emph{Management Science} 60(10):2583--2601.

\bibitem[{Goyal et~al.(2016)Goyal, Levi, \protect\BIBand{}
  Segev}]{goyal2016near}
Goyal V, Levi R, Segev D (2016) Near-optimal algorithms for the assortment
  planning problem under dynamic substitution and stochastic demand.
  \emph{Operations Research} 64(1):219--235.

\bibitem[{Gupte et~al.(2013)Gupte, Ahmed, Cheon, \protect\BIBand{}
  Dey}]{gupte2013solving}
Gupte A, Ahmed S, Cheon MS, Dey S (2013) Solving mixed integer bilinear
  problems using milp formulations. \emph{SIAM Journal on Optimization}
  23(2):721--744.

\bibitem[{Honhon et~al.(2010)Honhon, Gaur, \protect\BIBand{}
  Seshadri}]{honhon2010assortment}
Honhon D, Gaur V, Seshadri S (2010) Assortment planning and inventory decisions
  under stockout-based substitution. \emph{Operations research}
  58(5):1364--1379.

\bibitem[{Honhon et~al.(2012)Honhon, Jonnalagedda, \protect\BIBand{}
  Pan}]{honhon2012optimal}
Honhon D, Jonnalagedda S, Pan XA (2012) Optimal algorithms for assortment
  selection under ranking-based consumer choice models. \emph{Manufacturing \&
  Service Operations Management} 14(2):279--289.

\bibitem[{Jagabathula \protect\BIBand{}
  Rusmevichientong(2016)}]{jagabathula2016nonparametric}
Jagabathula S, Rusmevichientong P (2016) A nonparametric joint assortment and
  price choice model. \emph{Management Science} 63(9):3128--3145.

\bibitem[{K{\"o}k \protect\BIBand{} Fisher(2007)}]{kok2007demand}
K{\"o}k AG, Fisher ML (2007) Demand estimation and assortment optimization
  under substitution: Methodology and application. \emph{Operations Research}
  55(6):1001--1021.

\bibitem[{K{\"o}k et~al.(2008)K{\"o}k, Fisher, \protect\BIBand{}
  Vaidyanathan}]{kok2008assortment}
K{\"o}k AG, Fisher ML, Vaidyanathan R (2008) Assortment planning: Review of
  literature and industry practice. \emph{Retail supply chain management},
  99--153 (Springer).

\bibitem[{K{\"o}k et~al.(2015)K{\"o}k, Fisher, \protect\BIBand{}
  Vaidyanathan}]{kok2015assortment}
K{\"o}k AG, Fisher ML, Vaidyanathan R (2015) Assortment planning: Review of
  literature and industry practice. \emph{Retail supply chain management},
  175--236 (Springer).

\bibitem[{Kunnumkal(2015)}]{kunnumkal2015upper}
Kunnumkal S (2015) On upper bounds for assortment optimization under the
  mixture of multinomial logit models. \emph{Operations Research Letters}
  43(2):189--194.

\bibitem[{Kunnumkal \protect\BIBand{} Mart{\'\i}nez-de
  Alb{\'e}niz(2019)}]{kunnumkal2019tractable}
Kunnumkal S, Mart{\'\i}nez-de Alb{\'e}niz V (2019) Tractable approximations for
  assortment planning with product costs. \emph{Operations Research}
  67(2):436--452.

\bibitem[{Li \protect\BIBand{} Ke(2019)}]{li2019robust}
Li X, Ke J (2019) Robust assortment optimization using worst-case cvar under
  the multinomial logit model. \emph{Operations Research Letters}
  47(5):452--457.

\bibitem[{Liu et~al.(2019)Liu, Ma, \protect\BIBand{}
  Topaloglu}]{liu2019assortment}
Liu N, Ma Y, Topaloglu H (2019) Assortment optimization under the multinomial
  logit model with sequential offerings. \emph{INFORMS Journal on Computing (to
  appear)}
  \urlprefix\url{https://pdfs.semanticscholar.org/122c/4c3474d24b38c5baebd84d25e77265a728a4.pdf}.

\bibitem[{Ma et~al.(2019)Ma, Rusmevichientong, \protect\BIBand{}
  Topaloglu}]{ma2019assortment}
Ma Y, Rusmevichientong P, Topaloglu H (2019) Assortment optimization and
  pricing under the multinomial logit model with impatient customers. Technical
  report, Cornell University,
  \urlprefix\url{https://people.orie.cornell.edu/huseyin/publications/impatient_mnl.pdf}.

\bibitem[{Mahajan \protect\BIBand{} Van~Ryzin(2001)}]{mahajan2001stocking}
Mahajan S, Van~Ryzin G (2001) Stocking retail assortments under dynamic
  consumer substitution. \emph{Operations Research} 49(3):334--351.

\bibitem[{Perakis \protect\BIBand{} Roels(2010)}]{perakis2010robust}
Perakis G, Roels G (2010) Robust controls for network revenue management.
  \emph{Manufacturing \& Service Operations Management} 12(1):56--76.

\bibitem[{Rooderkerk \protect\BIBand{} van Heerde(2016)}]{rooderkerk2016robust}
Rooderkerk RP, van Heerde HJ (2016) Robust optimization of the 0--1 knapsack
  problem: Balancing risk and return in assortment optimization. \emph{European
  Journal of Operational Research} 250(3):842--854.

\bibitem[{Rusmevichientong et~al.(2010)Rusmevichientong, Shen,
  \protect\BIBand{} Shmoys}]{rusmevichientong2010dynamic}
Rusmevichientong P, Shen ZJM, Shmoys DB (2010) Dynamic assortment optimization
  with a multinomial logit choice model and capacity constraint.
  \emph{Operations research} 58(6):1666--1680.

\bibitem[{Rusmevichientong \protect\BIBand{}
  Topaloglu(2012)}]{rusmevichientong2012robust}
Rusmevichientong P, Topaloglu H (2012) Robust assortment optimization in
  revenue management under the multinomial logit choice model. \emph{Operations
  Research} 60(4):865--882.

\bibitem[{Rusmevichientong et~al.(2006)Rusmevichientong, Van~Roy,
  \protect\BIBand{} Glynn}]{rusmevichientong2006nonparametric}
Rusmevichientong P, Van~Roy B, Glynn PW (2006) A nonparametric approach to
  multiproduct pricing. \emph{Operations Research} 54(1):82--98.

\bibitem[{Ryzin \protect\BIBand{} Mahajan(1999)}]{ryzin1999relationship}
Ryzin Gv, Mahajan S (1999) On the relationship between inventory costs and
  variety benefits in retail assortments. \emph{Management Science}
  45(11):1496--1509.

\bibitem[{Smith \protect\BIBand{} Agrawal(2000)}]{smith2000management}
Smith SA, Agrawal N (2000) Management of multi-item retail inventory systems
  with demand substitution. \emph{Operations Research} 48(1):50--64.

\bibitem[{Talluri \protect\BIBand{} Van~Ryzin(2004)}]{talluri2004revenue}
Talluri K, Van~Ryzin G (2004) Revenue management under a general discrete
  choice model of consumer behavior. \emph{Management Science} 50(1):15--33.

\bibitem[{Vielma et~al.(2010)Vielma, Ahmed, \protect\BIBand{}
  Nemhauser}]{vielma2010mixed}
Vielma JP, Ahmed S, Nemhauser G (2010) Mixed-integer models for nonseparable
  piecewise-linear optimization: Unifying framework and extensions.
  \emph{Operations research} 58(2):303--315.

\bibitem[{Wang \protect\BIBand{} Sahin(2017)}]{wang2017impact}
Wang R, Sahin O (2017) The impact of consumer search cost on assortment
  planning and pricing. \emph{Management Science} 64(8):3649--3666.

\end{thebibliography}
\normalsize
\newpage
\pagestyle{empty}
\section*{Supplementary Material}
\vspace{3mm}

We present the updated version of Algorithms \ref{alg:2} and \ref{alg:greedy} considering multiple categories of customers in Algorithms \ref{alg:4} and \ref{alg:5}, respectively. 

\begin{algorithm}
	\small{
		\caption{: \PLUB{} approach for the constrained \APP{} with multiple customer categories  \label{alg:4}}
		\begin{algorithmic}[1]
			\State Set $ T =0 $, $ \theta_{min} = \theta(\emptyset) $, and $ \theta_{max} = 0 $
			\State Solve upper bound Problem (\ref{equ:equ51}) to obtain an initial optimal assortment $ S^\ast_T $ and $\mathbf{x}^\ast$, $ \boldsymbol{\alpha}^\ast$, and $ \boldsymbol{\beta}^\ast $ values
			\ForAll {$ c \in C $}
			\State Set $ \theta_c(S^\ast_T) = -\sum_{k=1}^{\bar u} \alpha^\ast_{kc}-\sum_{i =1}^n(1-x^\ast_i)\beta^\ast_{ic}  $ and $ \pi^\prime_{c} = h(\theta_c(S^\ast_T)) $
			\EndFor
			\While {$\pi ^ \prime _{c} - \phi (\theta_c(S^\ast_T)) > 0$ for some $ c \in C $} 
			\State Set $ T = T + 1 $ 
			\ForAll {$ c \in C $}
			\State Construct set of constraints $ Q_T^c $
			\EndFor
			\State Solve upper bound Problem $ (\ref{equ:equ51}) $ to obtain optimal assortment $ S^\ast_T $ and $\mathbf{x}^\ast$, $ \boldsymbol{\alpha}^\ast$, and $ \boldsymbol{\beta}^\ast $ values
			\ForAll {$ c \in C $}
			\State Set $ \theta_c(S^\ast_{T}) = -\sum_{k=1}^{\bar u} \alpha^\ast_{kc}-\sum_{i =1}^n(1-x^\ast_i)\beta^\ast_{ic}  $
			\EndFor
			\EndWhile
			\State \textbf{return} $ S_T^\ast$
	\end{algorithmic}}
\end{algorithm}
\normalsize  
\begin{algorithm} 
	\small{
		\caption{: Greedy algorithm for the constrained \APP{} with multiple customer categories} \label{alg:5}
		\begin{algorithmic}[1]
			\State Initialize $\underline z = 0$
			\ForAll{$i \in I \setminus \{0\}$}
			\State Initialize $S_i=\{i\}$ and $ \hat k = i $
			\ForAll {$ c \in C $}
			\State Use linear program $ (\ref{equ:equ24}) $ to find $ \theta_c(S_i) $ and set $ \pi_c(S_i) = e^ {\theta_c(S_i)} $ 
			\EndFor       
			\State Set $ \underline z= \sum_{c \in C} \bar w_c \left( \pi_c(S_i) \sum_{j \in S_i} \rho_{jc}(S_i)r_j \right) $
			\While {$ |S_i| \leq \bar c $ and $ \hat k \neq \texttt{NIL} $}
			\State Set $ \hat k = \texttt{NIL} $
			\ForAll{$k \in I \setminus (S_i \cup \{0\})$}
			\ForAll {$ c \in C $}
			\State Use linear program $ (\ref{equ:equ24}) $ to find $ \theta_c(S_i \cup \{k\}) $ and set $ \pi_c(S_i \cup \{k\}) = e^ {\theta_c(S_i \cup \{k\})} $
			\EndFor
			\State Compute $ z = \sum_{c \in C} \bar w_c \left(  \pi_c(S_i \cup \{k\}) \sum_{j \in S_i \cup \{k\}} \rho_{jc}(S_i \cup \{k\})r_j \right) $
			\If {$ z>\underline z $}
			\State Set $ \underline z=z $ and $ \hat k = k $
			\EndIf
			\EndFor
			\State Set $ S_i = S_i \cup \{\hat k\} $
			\EndWhile
			\EndFor
			\State Set $ i^\ast = \argmax_{i \in I \setminus \{0\}} \sum_{c \in C} \bar w_c \left( \pi_c(S_i) \sum_{j \in S_i} \rho_{jc}(S_i)r_j \right) $
			\State \textbf{return} $ S_{i^\ast} $
	\end{algorithmic}}
\end{algorithm}
\normalsize  

\end{document}